\documentclass[11pt]{extarticle}
\RequirePackage[OT1]{fontenc}
\RequirePackage{amsmath}
\usepackage{amssymb}
\RequirePackage{amsthm}
\usepackage{authblk}
\usepackage{caption}
\usepackage{comment}
\usepackage{dsfont}
\usepackage{eurosym}
\usepackage{geometry}
\usepackage{natbib}
\usepackage[colorlinks=true, allcolors=blue]{hyperref}
\usepackage{mathrsfs}
\usepackage{mathtools}
\usepackage{upquote}
\usepackage[dvipsnames]{xcolor}

\theoremstyle{plain}

\newtheorem{assumption}{Assumption}[section]

\newtheorem{proposition}{Proposition}[section]
\newtheorem{remark}{Remark}[section]
\newtheorem{theorem}{Theorem}[section]
\newtheorem{corollary}{Corollary}[section]
\theoremstyle{remark}

\newtheorem{example}{Example}[section]

\DeclareMathOperator{\sgn}{\mathrm{sgn}}

\DeclareMathOperator{\Erw}{\mathds E}

\DeclareMathOperator{\Prob}{\mathds P}

\usepackage{lipsum}
\usepackage{comment}

\begin{document}
\title{Estimating Lagged (Cross-)Covariance Operators of $L^p$-$m$-approximable Processes in Cartesian Product Hilbert Spaces}

\author[1]{Sebastian K\"uhnert\thanks{Sebastian.Kuehnert@ruhr-uni-bochum.de}}
\affil[1]{Fakultät für Mathematik, Ruhr-Universität Bochum, Bochum, DE}

\date{August 23, 2024}

\maketitle

\vspace{-.25in}

\begin{abstract} Estimating parameters of functional
ARMA, GARCH and invertible processes requires estimating lagged covariance and cross-covariance operators of Cartesian product Hilbert space-valued processes. Asymptotic results have been derived in recent years, either less generally or 
under a strict condition. This article derives upper bounds of the estimation errors for such operators based on the mild condition $L^p$-$m$-approximability for each lag, Cartesian power(s) and sample size, where the two processes can take values in different spaces in the context of lagged cross-covariance operators. Implications of our results on eigenelements, parameters in functional AR(MA) models and other general situations are also discussed. 
\end{abstract}



\vspace{-.05in}

\section{Introduction}\label{Section: Introduction}

Functional data analysis, which studies random functions or objects with a more general structure, has become increasingly important in recent decades, see \cite{HsingEubank2015},  \cite{HorvathKokoszka2012}, \cite{FerratyVieu2006} and \cite{RamsaySilverman2005} for extensive reviews. Examples of functional data are daily electricity prices and intra-day return curves, see \cite{Liebl2013} and \cite{Riceetal2020}, which are often assumed to have values in the Banach space $C[0,1],$ or the Hilbert space $L^2[0,1]$ of continuous and quadratic-Lebesgue integrable functions with domain $[0,1],$  respectively. When such data are observed continuously, e.g. \frenchspacing daily, they are called functional
time series (fTS). Excellent introductions to fTS analysis (fTSA) are \cite{Bosq2000} and \cite{kokoszka:2017:FDA-book}, and for useful concepts for common time series, see  \cite{ShumwayStoffer2017} and \cite{BrockwellDavis1991}. The dependence structure is one of the main objects of investigation in fTSA, measured by lagged covariance and cross-covariance operators, the functional counterparts of autocovariance and cross-covariance functions, respectively. Thus, there is a general interest in estimating these operators for which a reasonable (weak) dependence condition is required, see \cite{Dedecker2007} for a thorough review on weak dependence concepts. In fTSA, weak dependence in sense of cumulant mixing by  \cite{panaretos:tavakoli:2013} is frequently used, but $L^p$-$m$-approximability by \cite{HoermannKokoszka2010} even more often, as it usually follows from simple model assumptions. $L^p$-$m$-approximability is used e.g. \frenchspacing in functional linear regression, see \cite{Reimherr2015} and \cite{HoermannKidzinski2015},
and parameter estimation of common fTS in separable Hilbert spaces $\mathcal{H},$ such as functional \emph{autoregressive} (fAR), \emph{moving average} (fMA), ARMA (fARMA), (\emph{generalized}), \emph{autoregressive conditional heteroskedastic} (f(G)ARCH), and invertible linear processes, see \cite{Aueetal2017}, \cite{Bosq2000}, \cite{Hoermannetal2013}, \cite{Kuehnert2020}, \cite{KuehnertRiceAue2024+}, \cite{kuenzer:2024}, which requires the estimation of lagged covariance and cross-covariance operators of processes in Cartesian product Hilbert spaces. We would also like to point out the recently introduced and promising concept of \emph{$\pi$-dependence} by \cite{Kutta2024+}, which is a weak dependence concept for time series in general metric spaces with the property that it follows from both $\beta$-mixing and $L^p$-$m$-approximability.

Upper bounds of lagged covariance and cross-covariance operators' estimation errors  of strictly  stationary $\mathcal{H}$-valued processes have already been derived. \cite{HorvathKokoszka2012} established upper bounds for covariance operators of independent and identically distributed (i.i.d.) and \cite{HoermannKokoszka2010} for $L^4$-$m$-approximable processes in $L^2[0,1].$ Further, asymptotic upper bounds under $L^4$-$m$-approximibility were derived by \cite{AueKlepsch2017} for covariance operators in Cartesian products of processes originating from linear processes in $L^2[0,1],$  for covariance and cross-covariance operators by  \cite{HoermannKidzinski2015}, and for covariance and certain lagged covariance operators by \cite{kuenzer:2024}. Moreover, \cite{Kuehnert2022} established limits for the lagged covariance and cross-covariance operators' estimation errors for quite general processes in Cartesian products under a rather strict condition. Finally, for a comprehensive study on lagged cross-covariance operators of processes in $L^2[0,1],$ we refer to \cite{RiceShum2019}. 

This article establishes upper bounds  of estimation errors for \emph{lagged cross-covariance operators} $\mathscr{C}^h_{\!\boldsymbol{X}^{[m]},\boldsymbol{Y}^{[n]}}$ of the Cartesian space-valued processes $\boldsymbol{X}^{[m]}\coloneqq (X^{[m]}_k)_{k\in\mathbb{Z}}\subset\mathcal{H}^m$ and $\boldsymbol{Y}^{[n]}\coloneqq (Y^{[n]}_k)_{k\in\mathbb{Z}}\subset\mathcal{H}^n_\star$ in  \eqref{process X in Cartesian space}--\eqref{process Y in Cartesian space}, with $h\in\mathbb{Z},$ $m, n \in\mathbb{N},$ where the underlying processes  $\boldsymbol{X}\coloneqq(X_k)_{k\in\mathbb{Z}}$ and $\boldsymbol{Y}\coloneqq (Y_k)_{k\in\mathbb{Z}}$ are stationary and $L^4$-$m$-approximable with values in the real, separable Hilbert spaces $\mathcal{H}$ and $\mathcal{H}_\star$ with inner products $\langle\cdot, \cdot\rangle$ and $\langle\cdot, \cdot\rangle_\star,$ and induced norms  $\|\cdot\|$ and $\|\cdot\|_\star,$ respectively. All upper bounds are stated for any lag $h,$ Cartesian power(s) $m,n$ and sample size $N,$ where $h,m,n$ may depend on $N,$ and where no particular structure of $(X_k)$ and $(Y_k)$ is required. Further, for   $\boldsymbol{X}=\boldsymbol{Y}$ almost surely (a.s.) and for \emph{lagged covariance operators} $\mathscr{C}^h_{\!\boldsymbol{X}^{[m]}},$ the bounds can be slightly improved. We also discuss applications to the estimation of eigenelements, fAR(MA) parameters, and general sums and products involving lagged cross-covariance operators.

The remaining sections are organized as follows.   Section \ref{Section: Preliminaries} states our notation, basic concepts and $L^p$-$m$-approximability. Section \ref{Section: Main results} presents our main findings. Section \ref{Section: Consequences} states implications of our results. Section \ref{Conclusion} concludes the paper. Our proofs are collected in   Appendix \ref{Section Proofs}.

Throughout the whole article, convergences, $\mathrm o(\cdot)$ and $\mathrm O(\cdot)$ are for the sample size $N\rightarrow\infty,$ unless stated otherwise.

\vspace{-.05in}

\section{Preliminaries}\label{Section: Preliminaries}

\vspace{-.025in}

\subsection{Notation and basic definitions}

First we write $\mathbb{N}_0\coloneqq \mathbb{N}\cup\{0\}.$ For the \emph{indicator function} of a subset $A\subset\mathbb{R},$ we write $\mathds 1_A(x),$ $ x\in\mathbb{R}.$ The additive identity
of linear spaces is denoted by $0,$ and the identity map by  $\mathbb{I}.$ The Cartesian product space $\mathcal{H}^n,$ with component-wise scalar multiplication and addition, is a Hilbert space with inner product $\langle \boldsymbol{x} ,\boldsymbol{y}\rangle = \sum_{i=1}^n\langle x_i,y_i\rangle$ and norm $\|\boldsymbol{x}\| = (\sum^n_{i=1}\|x_i\|^2)^{1/2},$ with $\boldsymbol{x} \coloneqq (x_1, \dots, x_n)^{\!\top}\!, \boldsymbol{y} \coloneqq (y_1, \dots, y_n)^{\!\top} \in \mathcal{H}^n$. Further, $\mathcal{L}_{\mathcal{H}, \mathcal{H}_\star}, \mathcal{S}_{\mathcal{H}, \mathcal{H}_\star}\subsetneq\mathcal{L}_{\mathcal{H}, \mathcal{H}_\star},$  and $\mathcal{N}_{\mathcal{H}, \mathcal{H}_\star}\subsetneq\mathcal{S}_{\mathcal{H}, \mathcal{H}_\star}$ denote the spaces of linear, bounded; Hilbert-Schmidt (H-S); and nuclear (trace-class) operators $A\colon \mathcal{H} \rightarrow \mathcal{H}_\star,$ respectively, with operator norm  $\|\cdot\|_{\mathcal{L}};$ H-S inner product $\langle\cdot,\cdot \rangle_{\mathcal{S}}$ and H-S norm $\|\cdot\|_{\mathcal{S}};$ and nuclear norm $\|\cdot\|_{\mathcal{N}},$  respectively, with $\mathcal{T_H}\coloneqq\mathcal{T}_{\mathcal{H, H}}$ for $\mathcal{T}\in\{\mathcal{L,S,N}\}.$ All our random variables are defined on a common probability space $(\Omega, \mathfrak{A},\mathbb{P}).$ We write $X\stackrel{d}{=} Y$ for identically distributed random variables $X$ and $Y.$ For $p \in [1,\infty),$ $L^p_{\mathcal{H}} = L^p_{\mathcal{H}}(\Omega, \mathfrak{A},\mathbb{P})$ is the space of (classes of) random variables $X \in \mathcal{H}$ with $\nu_p(X) \coloneqq (\Erw\!\|X\|^p)^{1/p} < \infty,$  and a process $(X_{\hspace{-0.05em}k})\subset L^1_{\mathcal{H}}$ is \emph{centered}, if $\Erw(X_{\hspace{-0.05em}k}) = 0$ for all $k,$ with expectation in the Bochner-integral sense. The \emph{cross-covariance operator} of $X\in L^2_{\mathcal{H}}$ and $Y \in L^2_{\mathcal{H}_\star}$ is defined by 
$\mathscr{C}_{\!X,Y}\coloneqq\Erw[X-\Erw(X)]\!\otimes\![Y-\Erw(Y)],$ and the \emph{covariance operator} of $X$ by $\mathscr{C}_{\!X} = \mathscr{C}_{\!X,X},$ where $x\otimes y\coloneqq \langle x, \cdot\rangle y\in\mathcal{S}_{\mathcal{H},\mathcal{H}_\star}$ is the \emph{tensorial product} of $x\in\mathcal{H}, y\in\mathcal{H}_\star.$  Further, $\mathscr{C}_{\!X,Y}\! \in \mathcal{N}_{\mathcal{H},\mathcal{H}_\star}$ and $\mathscr{C}^\ast_{\!X,Y}\! =\mathscr{C}_{Y,X}\in\mathcal{N}_{\mathcal{H}_\star,\mathcal{H}},$ thus $\mathscr{C}_{\!X}\! \in \mathcal{N_H}$ and $\mathscr{C}^\ast_{\!X}\! =\mathscr{C}_{X},$ where $A^{\ast}$ is the adjoint of $A\in\mathcal{L}_{\mathcal{H}, \mathcal{H}_\star}.$ 

A process $\boldsymbol{X}\!=\!(X_k)_{k\in\mathbb{Z}} \subset \mathcal{H}$ is  \emph{strictly stationary} if $(X_{t_1}, \dots, X_{t_n})\stackrel{d}{=} (X_{t_1 + h}, \dots, X_{t_n +h})$ for any $h, t_1, \dots, t_n\in\mathbb{Z}, n\in\mathbb{N},$ and \emph{weakly stationary} if $\boldsymbol{X}\!\subset\! L^2_{\mathcal{H}},$ $\Erw(X_k) = \Erw(X_\ell)$ for any $k,\ell\in\mathbb{Z},$ and $\mathscr{C}_{X_k, X_\ell} = \mathscr{C}_{X_{k+h}, X_{\ell+h}}$ for any $h,k,\ell\in\mathbb{Z}.$ For $h\in \mathbb{Z},$ the \emph{lag-$h$-covariance operators} of a weakly stationary process $\boldsymbol{X}\!\subset\! L^2_{\mathcal{H}}$ are defined by $\mathscr{C}^h_{\!\boldsymbol{X}} \coloneqq \mathscr{C}_{\!X_0,X_h},$ and we call $\mathscr{C}_{\!\boldsymbol{X}}\coloneqq\mathscr{C}^0_{\!\boldsymbol{X}}$ covariance operator of $\boldsymbol{X}.$ Two processes $(X_k)_{k\in\mathbb{Z}}\!\subset\!\mathcal{H}, (Y_k)_{k\in\mathbb{Z}}\!\subset\!\mathcal{H}_\star$ are \emph{jointly strictly stationary} if $(X_{s_1}, \dots, X_{s_m}, Y_{t_1}, \dots, Y_{t_n})\stackrel{d}{=} (X_{s_1+h}, \dots, X_{s_m+h}, Y_{t_1+h}, \dots, Y_{t_n+h})$ for all $h, s_1, \dots, s_m, t_1, \dots, t_n\in\mathbb{Z}, m,n\in\mathbb{N}.$ %This holds, e.g., if $(X_k)_k, (Y_k)_k$ are i.i.d. \frenchspacing processes that are independent from each other. 
Further, two processes $\boldsymbol{X}=(X_k)_{k\in\mathbb{Z}}\subset L^2_{\mathcal{H}},$ $ \boldsymbol{Y}=(Y_k)_{k\in\mathbb{Z}}\subset L^2_{\mathcal{H}_\star}$ are \emph{jointly weakly stationary} if $\boldsymbol{X}$ and $\boldsymbol{Y}$ are weakly stationary and $\mathscr{C}_{X_k, Y_\ell} = \mathscr{C}_{X_{k+h}, Y_{\ell+h}}$ for any  $h,k,\ell,$ and the  \emph{lag-$h$-cross covariance operators} of such processes  $\boldsymbol{X}, \boldsymbol{Y}$ are defined by $\mathscr{C}^h_{\!\boldsymbol{X}\!,\boldsymbol{Y}}\coloneqq 	\mathscr{C}_{\!X_0,Y_h}, h\in\mathbb{Z}.$ % Joint weak stationarity holds, e.g., for two fAR$(1)$ processes defined by $X_k=\gamma(X_{k-1}) + \varepsilon_k$ a.s.,    $Y_k=\delta(Y_{k-1}) + \varepsilon_k$ a.s. \frenchspacing  for each $k,$ with $\|\gamma\|_{\mathcal{L}}, \|\delta\|_{\mathcal{L}} <1,$ and where $(\varepsilon_k)_k\subset L^2_\mathcal{H}$ is i.i.d. \frenchspacing and centered.

\vspace{-.025in}

\subsection{Comments on $L^p$-$m$-approximability}\label{Comments Lp-m}
A process $(X_k)_{k\in\mathbb{Z}}\subset\mathcal{H}$ is called  \emph{$L^p_{\hspace{-0.05em}\mathcal{H}}$-$m$-approximable} for $p\geq 1$ if for each $k$ holds $X_k\in L^p_{\mathcal{H}}$ and $X_k = f(\varepsilon_k, \varepsilon_{k-1}, ...),$ where $(\varepsilon_k)_{k\in\mathbb{Z}}$ is an i.i.d. \frenchspacing process in a measurable space $S,$ and $f\colon S^\infty \rightarrow \mathcal{H}$ is a measurable function, i.e. \frenchspacing $(X_k)$ is \emph{causal} (with regard  to (w.r.t) $(\varepsilon_k)$), such that
\begin{align*}
\sum_{m=1}^\infty\, 
\nu_p(X_{\hspace{-0.05em}m}- X^{(m)}_{ 
	\hspace{-0.05em}m})\, < \infty,
\end{align*} 
with  $X^{(\ell)}_k\coloneqq f(\varepsilon_k, \dots, \varepsilon_{k-\ell+1}, \varepsilon^{(\ell)}_{k-\ell}, \varepsilon^{(\ell)}_{k-\ell-1}, \dots)$ for $k\in\mathbb{Z},$ $\ell\in\mathbb{N}_0,$ and independent copies $(\varepsilon^{(\ell)}_k)_k$ of $(\varepsilon_k).$ Thereby, $L^p$-$m$-approximability implies stationarity, see \cite{Stout1974}, and  $(X^{(m)}_k)_k$ are $m$-dependent for each $m,$ with $X^{(m)}_k\stackrel{d}{=}X_k$ for each $k,m.$ That $L^p$-$m$-approximability is feasible can be observed, e.g., on the basis of the general class of causal, \emph{linear processes} $X_k=\sum_{i\in\mathbb{N}}\phi_i(\varepsilon_{k-i}) + \varepsilon_k$ a.s. \frenchspacing for each $k,$  with $\phi_i\in\mathcal{L_H}$ and i.i.d., centered process $(\varepsilon_k)_{k\in\mathbb{Z}}\subset L^2_{\mathcal{H}},$ where the series defining $X_k$ converges a.s. \frenchspacing and in $L^2_{\mathcal{H}}$ if  $\sum_{i\in\mathbb{N}}\|\phi_i\|^2_{\mathcal{L}} < \infty,$ see \cite{Bosq2000}. For $p\geq 2,$ $L^p$-$m$-approximability of a linear process $(X_k)\subset L^2_{\mathcal{H}}$ holds according to \cite{HoermannKokoszka2010}, Proposition 2.1, if  $\nu_p(\varepsilon_0) = (\Erw\!\|\varepsilon_0\|^p)^{1/p} < \infty$ and $\sum_{m\in\mathbb{N}}\sum_{j\geq m}\|\phi_j\|_{\mathcal{L}} < \infty.$ For an fAR$(1)$ process 
\begin{align}
  X_k = \psi(X_{k-1}) + \varepsilon_k\;\,\text{a.s.},\quad k\in\mathbb{Z},\label{equation fAR(1)}   
\end{align}
with $\psi\in\mathcal{L_H}$ and i.i.d., centered process $(\varepsilon_k)\subset L^2_{\mathcal{H}},$ $L^p$-$m$-approximability holds under a quite simple condition. 

\begin{proposition}\label{Prop. L^p-m for AR} Let $(X_k)$ be the fAR$(1)$ process in \eqref{equation fAR(1)}. Further, suppose   $\|\psi^{j_0}\|_{\mathcal{L}}<1$ for some $j_0\in\mathbb{N},$ and  $\nu_p(\varepsilon_0) <\infty$ for some $p\geq 2.$ Then,  $(X_k)$ is $L^p$-$m$-approximable.
\end{proposition}

\begin{example}\label{Ex: upper bound for L^p errors fAR(1)} Let $(X_k)$ be the fAR$(1)$ process  in \eqref{equation fAR(1)}, with $\xi\coloneqq\|\psi\|_{\mathcal{L}}<1$ and $\nu_p(\varepsilon_0) <\infty$ for some $p\geq2.$ Then,   $\|\psi^j\|_{\mathcal{L}}\leq   \xi^j$ for $j\in\mathbb{N},$ and $0\leq \xi < 1$ yield (see the proof of Proposition \ref{Prop. L^p-m for AR}):
\begin{align*}
\sum^\infty_{m=1}\nu_p\big(X_m - X^{(m)}_m\big) \leq \frac{2\nu_p(\varepsilon_0)}{1-\xi}\sum^\infty_{m=1}\xi^m = \frac{2\xi}{(1-\xi)^2}\,\nu_p(\varepsilon_0) < \infty.
\end{align*}
\end{example}

\vspace{-.05in}

\section{Main results}\label{Section: Main results}

We want to derive  estimates for the lag-$h$-cross-covariance operators $\mathscr{C}^h_{\!\boldsymbol{X}^{[m]}, \boldsymbol{Y}^{[n]}}\!\in\!\mathcal{N}_{\mathcal{H}^m,\mathcal{H}^n}$ and subsequently for the lag-$h$-covariance operators $\mathscr{C}^h_{\!\boldsymbol{X}^{[m]}}\!\in\!\mathcal{N}_{\mathcal{H}^m},$ with $h\in\mathbb{Z}$ and $m,n\in\mathbb{N},$ where  $\boldsymbol{X}^{[m]}=(X^{[m]}_k)_{k\in\mathbb{Z}}\subset\mathcal{H}^m$ and $\boldsymbol{Y}^{[n]}=(Y^{[n]}_k)_{k\in\mathbb{Z}}\subset\mathcal{H}^n_\star$ are the Cartesian product Hilbert space-valued processes whose elements are defined by
\begin{alignat}{2}
  X^{[m]}_k&\coloneqq (X_k, X_{k-1}, \dots, X_{k-m+1})^\top\!, &&\quad k\in\mathbb{Z},\label{process X in Cartesian space}\\
  Y^{[n]}_k&\coloneqq (Y_k, Y_{k-1}, \dots, Y_{k-n+1})^\top\!, &&\quad k\in\mathbb{Z}.\label{process Y in Cartesian space}
\end{alignat}
Estimating the lag-$h$-cross-covariance operators of these processes involves the approximates  
\begin{alignat}{2}
  X^{[m],(\ell)}_k&\coloneqq (X^{(\ell)}_k, X^{(\ell)}_{k-1}, \dots, X^{(\ell)}_{k-m+1})^\top\!, &&\quad k\in\mathbb{Z},\,\ell\in\mathbb{N}_0,\label{approximate of process X in Cartesian space}\\
  Y^{[n],(\ell)}_k&\coloneqq (Y^{(\ell)}_k, Y^{(\ell)}_{k-1}, \dots, Y^{(\ell)}_{k-n+1})^\top\!, &&\quad k\in\mathbb{Z},\,\ell\in\mathbb{N}_0,\label{approximate of process Y in Cartesian space}
\end{alignat}
provided $\boldsymbol{X}=(X_k)_{k\in\mathbb{Z}}\subset \mathcal{H}$ and $\boldsymbol{Y}=(Y_k)_{k\in\mathbb{Z}}\subset \mathcal{H}_\star$ satisfy the following.

\begin{assumption}\label{As: X L^2, centered} The processes  $\boldsymbol{X}$ and $\boldsymbol{Y}$ are   $L^4$-$m$-approximable and centered. Thereby, 
\begin{alignat*}{2}
  X_k&= f(\varepsilon_k, \varepsilon_{k-1}, \dots ), &&\quad k\in\mathbb{Z},\\
  Y_k&= g(\eta_k, \eta_{k-1}, \dots ), &&\quad k\in\mathbb{Z},
\end{alignat*}
where $(\varepsilon_k)_{k\in\mathbb{Z}}$ and $(\eta_k)_{k\in\mathbb{Z}}$ are i.i.d. \frenchspacing processes with values in the measurable space $S,$ where $\varepsilon_k$ and $\eta_\ell$ are independent for $k\neq \ell,$ and where $f\colon S^\infty \rightarrow \mathcal{H}$ and $g\colon S^\infty \rightarrow \mathcal{H}_\star$ are measurable functions.
\end{assumption}

\noindent To establish our estimates, joint weak stationarity, and for technical reasons, also joint weak stationarity of the fourth moments is required.

\begin{assumption}\label{As: jointly stationary} For each  $h,k,\ell$ holds $\mathscr{C}_{X_k, Y_\ell} = \mathscr{C}_{X_{k+h}, Y_{\ell+h}},$ and for each $h,i,j,k,\ell$ holds    
$$    
  \Erw\langle X_i\otimes Y_j, X_k\otimes Y_{\ell}\rangle_{\mathcal{S}} = \Erw\langle X_{i+h}\otimes Y_{j+h}, X_{k+h}\otimes Y_{\ell+h}\rangle_{\mathcal{S}}\,.
$$  
\end{assumption}

\begin{proposition}\label{Prop: Joint stat cart} Let Assumptions \ref{As: X L^2, centered}--\ref{As: jointly stationary} hold. Then, $\boldsymbol{X}^{[m]}$ and $\boldsymbol{Y}^{[n]}$ are jointly weakly stationary, and for each  $h,i,j,k,\ell,m, n$ holds 
$$\Erw\langle X^{[m]}_i\otimes Y^{[n]}_j, X^{[m]}_k\otimes Y^{[n]}_\ell\rangle_{\mathcal{S}} = \Erw\langle X^{[m]}_{i+h}\otimes Y^{[n]}_{j+h}, X^{[m]}_{k+h}\otimes Y^{[n]}_{\ell+h}\rangle_{\mathcal{S}}\,.$$ 
\end{proposition} 

\begin{remark} Assumption \ref{As: X L^2, centered} is quite mild, see Example \ref{Ex: upper bound for L^p errors fAR(1)}, as is Assumption  \ref{As: jointly stationary} which follows from joint strict stationarity of $\boldsymbol{X}\subset L^4_{\mathcal{H}}$ and $ \boldsymbol{Y}\subset L^4_{\mathcal{H}_\star}.$ Also, if $\boldsymbol{X}=\boldsymbol{Y}$ a.s., Assumption \ref{As: jointly stationary} follows from strict stationarity of $\boldsymbol{X},$ which is given by  Assumption \ref{As: X L^2, centered}.     
\end{remark}

\vspace{-.025in}

\subsection{Estimation of lagged cross-covariance operators}

Based on samples $X_1, \dots, X_N$ of $\boldsymbol{X}$ and $Y_1, \dots, Y_N$ of $\boldsymbol{Y}$ with sample size $N \in\mathbb{N},$ for each $h, m, n$ with $m,n\leq N$ and $n-N \leq h \leq N-m,$ and $N'\coloneqq\min\{N,N-h\} - \max\{m,n-h\}+1,$
\begin{align*}
  \hat{\mathscr{C}}^h_{\!\boldsymbol{X}^{[m]}, \boldsymbol{Y}^{[n]}}\coloneqq \frac{1}{N'}\sum^{\min\{N,N-h\}}_{k=\max\{m,n-h\}}X^{[m]}_k\otimes Y^{[n]}_{k+h}
\end{align*}
 is an unbiased estimator for $\mathscr{C}^h_{\!\boldsymbol{X}^{[m]}, \boldsymbol{Y}^{[n]}},$ with    $(\hat{\mathscr{C}}^h_{\!\boldsymbol{X}^{[m]}, \boldsymbol{Y}^{[n]}})^\ast = \hat{\mathscr{C}}^{-h}_{\!\boldsymbol{Y}^{[n]}, \boldsymbol{X}^{[m]}}.$

\begin{theorem}\label{Theo cross-cov op} Let Assumptions \ref{As: X L^2, centered}--\ref{As: jointly stationary} hold. Then, for each $h\in\mathbb{Z},$ $m, n, N\in\mathbb{N}$ with $m,n \leq N$ and $n-N \leq h\leq N-m,$ and with $\kappa\coloneqq \max\{m,n-h\} + \mathds 1_{\mathbb{N}}(h) \cdot h - 1$ and $\kappa'\coloneqq \max\{1,\kappa\}$ holds     
\begin{align}\label{claim cross-cov operator}
  \frac{N'}{mn(2\kappa'-1)}\Erw\!\big\|\hat{\mathscr{C}}^h_{\!\boldsymbol{X}^{[m]}, \boldsymbol{Y}^{[n]}}\! - \mathscr{C}^h_{\!\boldsymbol{X}^{[m]}, \boldsymbol{Y}^{[n]}}\big\|^2_{\mathcal{S}} \leq \,\xi_{\boldsymbol{X,Y}}(h,m,n),   
\end{align}
where $0 \leq \xi_{\boldsymbol{X,Y}}(h,m,n) < \infty$ depends for any $h, m,n$ on the fourth moments of $\boldsymbol{X}, \boldsymbol{Y}\!,$ and the $L^4$-$m$-approximation errors through 
\begin{align*}
  \xi_{\boldsymbol{X,Y}}(h,m,n) &\coloneqq \nu^2_4(X_1)\nu^2_4(Y_1)\notag\allowdisplaybreaks\\
  &\qquad + \frac{2\sqrt{2}}{2\kappa' - 1}\,\nu_4(X_1)\nu_4(Y_1)\sum_{k\geq \kappa'}\nu_4(Y_1)\nu_4\big(X_k - X^{(k)}_k\big) + \nu_4(X_1)\nu_4\big(Y_k - Y^{(k)}_k\big).
\end{align*}
\end{theorem} 
\noindent Theorem \ref{Theo cross-cov op} entails for $h=h_N,$ $m=m_N,$ $ n=n_N=\mathrm{o}(N)$ with $\max\{|h|, m, n\}\rightarrow\infty\colon$ 
\begin{align*}
  \limsup_{N\rightarrow\infty}\frac{N}{mn(2\kappa'-1)}\Erw\!\big\|\hat{\mathscr{C}}^h_{\!\boldsymbol{X}^{[m]}, \boldsymbol{Y}^{[n]}}\! - \mathscr{C}^h_{\!\boldsymbol{X}^{[m]}, \boldsymbol{Y}^{[n]}}\big\|^2_{\mathcal{S}} \leq \nu^2_4(X_1)\nu^2_4(Y_1)\,.
\end{align*}
Also, both for $\max\{|h|, m, n\}\rightarrow\infty$ and $\max\{|h|, m, n\}\not\rightarrow\infty,$ we have
\begin{align*}
  \Erw\!\big\|\hat{\mathscr{C}}^h_{\!\boldsymbol{X}^{[m]}, \boldsymbol{Y}^{[n]}}\! - \mathscr{C}^h_{\!\boldsymbol{X}^{[m]}, \boldsymbol{Y}^{[n]}}\big\|^2_{\mathcal{S}} = \mathrm O_{\Prob}(mn\max\{|h|,m,n\}N^{-1})\,.%\begin{cases}
    %\,\mathrm O_{\Prob}(N^{-1}), &\textnormal{if}\;h, m, n\;\text{are fixed},\\
    %\,\mathrm O_{\Prob}(hN^{-1}), &\textnormal{if}\;h\rightarrow\pm\infty\;\text{and}\;m, n\;\text{are fixed},\\
% \,\mathrm O_{\Prob}(m^2N^{-1}), &\textnormal{if}\;h\;\text{is fixed,}\;m\rightarrow\infty\;\text{and}\;n\;\text{is fixed},\\
% \,\mathrm O_{\Prob}(n^2N^{-1}), &\textnormal{if}\;h,m\;\text{are fixed and}\;n\rightarrow\infty,\\
% \,\mathrm O_{\Prob}(m\max\{|h|,m\}N^{-1}), &\textnormal{if}\;h\rightarrow\pm\infty, m\rightarrow\infty\;\text{and}\;n\;\text{is fixed},\\
% \,\mathrm O_{\Prob}(n\max\{|h|,n\}N^{-1}), &\textnormal{if}\;h\rightarrow\pm\infty,\,m\;\text{is fixed,}\,n\rightarrow\infty,\\
% \,\mathrm O_{\Prob}(mn\max\{m,n\}N^{-1}), &\textnormal{if}\;h\;\text{is fixed and}\;m,n\rightarrow\infty,\\
% \,\mathrm O_{\Prob}(mn\max\{|h|,m,n\}N^{-1}), &\textnormal{if}\;h\rightarrow\pm\infty\;\text{and}\;m, n\rightarrow\infty.
%  \end{cases} 
\end{align*}

Furthermore, if $\boldsymbol{X}=\boldsymbol{Y}$ a.s., our upper bound becomes 
\begin{align}\label{upper bound X=Y}
\xi_{\boldsymbol{X,X}}(h,m,n) &= \nu^4_4(X_1) + \frac{4\sqrt{2}}{2\kappa' - 1}\,\nu^3_4(X_1)\!\sum_{k\geq \kappa'}\nu_4\big(X_k - X^{(k)}_k\big)\,.
\end{align}
For $h=0,$ $m=n=1,$ i.e. \frenchspacing for covariance operators of $\mathcal{H}$-valued processes, this bound coincides with the 'classical' upper bound derived for $L^2[0,1]$-valued processes in \cite{HoermannKokoszka2010}. In fact, if $\boldsymbol{X}=\boldsymbol{Y}$ a.s., we actually have the following slightly  sharper bound. 
\begin{proposition}\label{Corollary Cross-cov on X's} Under Assumption \ref{As: X L^2, centered} holds for each $h, m, n, N$ with $m,n\leq N$ and $n-N \leq h\leq N-m,$ and with $\kappa\coloneqq \max\{m,n-h\} + \mathds 1_{\mathbb{N}}(h) \cdot h - 1$ and $\kappa'=\max\{1,\kappa\}\colon$      
\begin{align*}
    \frac{N'}{mn(2\kappa'-1)}\Erw\!\big\|\hat{\mathscr{C}}^h_{\!\boldsymbol{X}^{[m]}, \boldsymbol{X}^{[n]}}\! - \mathscr{C}^h_{\!\boldsymbol{X}^{[m]}, \boldsymbol{X}^{[n]}}\big\|^2_{\mathcal{S}} \leq \tau_{\boldsymbol{X}}(h,m,n),  
\end{align*}
where $0 \leq \tau_{\boldsymbol{X}}(h,m,n) < \infty$ depends for any $h, m,n$ on the fourth moments of $\boldsymbol{X},$ and the $L^4$-$m$-approximation errors through 
\begin{align}\label{Def tau_X}
\tau_{\boldsymbol{X}}(h,m,n) &= \nu^4_4(X_1) + \frac{4}{2\kappa' - 1}\,\nu^3_4(X_1)\!\sum_{k\geq \kappa'}\nu_4\big(X_k - X^{(k)}_k\big)\,.
\end{align}
\end{proposition}

\begin{remark}\label{Remark lag-h-cross-cov-op estimation}
\begin{enumerate}
    \item[\textnormal{(a)}] The upper bounds of the lagged cross-covariance operators'  estimation errors of Cartesian product Hilbert space-valued processes under $L^4$-$m$-approximability in Theorem \ref{Theo cross-cov op}, stated for any lag $h,$ Cartesian power $m,n$ and sample size $N,$ are, as far as we know, new in all their generality. However, \cite{HoermannKidzinski2015} already gave a result for certain lagged cross-covariance operators under $L^4$-$m$-approximability but without considering Cartesian products, and \cite{Kuehnert2022} deduced limits for arbitrary lagged cross-covariance operators in Cartesian products under the rather strict summability condition 
\begin{align}\label{summ cond.}
\sum_{i\in\mathbb{Z}}\sum_{j\in\mathbb{Z}}\,\big|\Erw\!\big\langle X_0\!\otimes\!Y_i - \mathscr{C}^{i}_{\!\boldsymbol{X,Y}}, X_j\!\otimes\!Y_{i+j} - \mathscr{C}^{j}_{\!\boldsymbol{X,Y}} \big\rangle_{\!\mathcal{S}}\big| < \infty.
\end{align} \item[\textnormal{(b)}] Instead of the approximates \eqref{approximate of process X in Cartesian space}--\eqref{approximate of process Y in Cartesian space} one could also use $X^{[m],(\ell)}_k\coloneqq (X^{(\ell)}_k, X^{(\ell-1)}_{k-1}, \dots, X^{(\ell-m+1)}_{k-m+1})^\top$ and $Y^{[n],(\ell)}_k\coloneqq (Y^{(\ell)}_k, Y^{(\ell-1)}_{k-1}, \dots, Y^{(\ell-n+1)}_{k-n+1})^\top\!.$ At a first glance using these approximates seems more advantageous, as independence of $Z^{[m,n]}_{1,h}=X^{[m]}_1\otimes Y^{[n]}_{1+h}$ and $Z^{[m,n],(\ell)}_{1+k,h}=X^{[m],(\ell)}_{1+k}\otimes Y^{[n],(\ell)}_{1+k+h}$ is given for $k>|h|,$ leading to improved upper bounds for fixed $h$ when $m\rightarrow\infty$ or $n\rightarrow\infty.$  Nevertheless, well-definedness of $Z^{[m,n],(\ell)}_{1+k,h}$ and the aforementioned independence is given for  any $\ell$ with $\max\{m,n\}-1\leq \ell \leq k-|h|,$ thus for $k\geq \max\{1,\max\{m,n\}+|h|-1\},$ hence, since all further transformations with these approximates would be the same as with the approximates    \eqref{approximate of process X in Cartesian space}--\eqref{approximate of process Y in Cartesian space} in the proof of Theorem \ref{Theo cross-cov op}, we would get comparable rates. \item[\textnormal{(c)}] The series expression in $\tau_{\boldsymbol{X}}(h,m,n)$ does not contain $\sqrt{2}$, contrary to $\xi_{\boldsymbol{X,X}}(h,m,n)$, because the Minkowski inequality is not needed in the proof of Proposition \ref{Corollary Cross-cov on X's}.
    \item[\textnormal{(d)}] The processes in \cite{Kuehnert2022} had a more general structure than  \eqref{process X in Cartesian space}--\eqref{process Y in Cartesian space}, and the sample sizes of the samples from  $\boldsymbol{X}$ and $\boldsymbol{Y}$ were allowed to be  different. For the sake of clarity, however, we will refrain from making further generalizations.
\end{enumerate}
\end{remark}

\vspace{-.025in}

\subsection{Estimation of lagged covariance operators}
Proposition \ref{Corollary Cross-cov on X's}  immediately leads to upper bounds for the estimation errors for lagged covariance operators by setting $m=n.$ Given a sample $X_1, \dots, X_N$ of $\boldsymbol{X},$ the estimator for the lag-$h$-covariance operator   $\mathscr{C}^h_{\!\boldsymbol{X}^{[m]}}$ for any $h, m$ with $m\leq N$ and $|h| \leq N-m$ has the simpler form 
\begin{align*}
  \hat{\mathscr{C}}^h_{\!\boldsymbol{X}^{[m]}, \boldsymbol{X}^{[m]}} = \hat{\mathscr{C}}^h_{\!\boldsymbol{X}^{[m]}} = 
    \frac{1}{N-|h|-m+1}\,\sum^{\min\{N,N-h\}}_{k=\max\{m,m-h\}}X^{[m]}_k\otimes X^{[m]}_{k+h}\,.
\end{align*}

\begin{corollary}\label{Cor cov op} Under Assumption \ref{As: X L^2, centered} holds for each  $h, m, N$ with $|h|\leq N-m\colon$ 
\begin{align*}
\frac{N\!-\!|h|\!-\!m\!+\!1}{m^2(2\kappa'-1)}\Erw\!\big\|\hat{\mathscr{C}}^h_{\!\boldsymbol{X}^{[m]}}\! - \mathscr{C}^h_{\!\boldsymbol{X}^{[m]}}\big\|^2_{\mathcal{S}} \leq \tilde{\tau}_{\boldsymbol{X}}(h,m), 
\end{align*}
with $\kappa'\coloneqq\max\{1,m+|h|-1\}$ and  $0 \leq \tilde{\tau}_{\boldsymbol{X}}(h,m)\coloneqq\tau_{\boldsymbol{X}}(h,m,m) < \infty$ defined in \eqref{Def tau_X}.  
\end{corollary}

\noindent From Corollary \ref{Cor cov op} follows for $h=h_N,$ $m=m_N=\mathrm o(N)$ with $\max\{|h|,m\}\rightarrow \infty\colon$  
\begin{alignat*}{2}
    \limsup_{N\rightarrow\infty}\frac{N}{2m^2(|h|+m)}\Erw\!\big\|\hat{\mathscr{C}}^h_{\!\boldsymbol{X}^{[m]}}\! - \mathscr{C}^h_{\!\boldsymbol{X}^{[m]}}\big\|^2_{\mathcal{S}} &\leq \nu^4_4(X_1)\,.
\end{alignat*}
Further, both for $\max\{|h|,m\}\rightarrow \infty$ and $\max\{|h|,m\}\not\rightarrow \infty,$ it holds 
\begin{align*} 
  \Erw\!\big\|\hat{\mathscr{C}}^h_{\!\boldsymbol{X}^{[m]}}\! - \mathscr{C}^h_{\!\boldsymbol{X}^{[m]}}\big\|^2_{\mathcal{S}} = \mathrm O_{\Prob}(m^2\!\max\{|h|,m\}N^{-1})\,.
  %\begin{cases}
  %    \,\mathrm O_{\Prob}(N^{-1}), &\textnormal{if}\;h, m\;\text{are fixed},\\
  %   \,\mathrm O_{\Prob}(hN^{-1}), &\textnormal{if}\;h\rightarrow\pm\infty\;\text{and}\;m\;\text{is fixed},\\
% \,\mathrm O_{\Prob}(m^3N^{-1}), &\textnormal{if}\;h\;\text{is fixed and}\;m\rightarrow\infty,\\
% \,\mathrm O_{\Prob}(m^2\!\max\{|h|,m\}N^{-1}), &\textnormal{if}\;h\rightarrow \pm\infty\;\text{and}\; m\rightarrow\infty.
%  \end{cases} 
\end{align*}
\begin{remark} As Corrolary \ref{Cor cov op} states upper bounds of lagged covariance operators'  estimation errors in Cartesian products for any lag $h,$ Cartesian power $m$ and sample size $N$ under $L^4$-$m$-approximability, it is new in several ways. Previously,  \cite{HoermannKokoszka2010} stated the somewhat cruder upper bound $\xi_{\boldsymbol{X,X}}(0,1,1)$ in \eqref{upper bound X=Y} instead of $\tilde{\tau}_{\boldsymbol{X}}(0,1)$  for covariance operators of $L^2[0,1]$-valued processes under $L^4$-$m$-approximability, see Remark \ref{Remark lag-h-cross-cov-op estimation}(c). Further,   \cite{AueKlepsch2017} derived upper bounds for covariance operators of $L^2[0,1]$-valued linear processes in Cartesian products, and \cite{Kuehnert2022} deduced limits for the lagged covariance operators' estimation errors in Cartesian products under a  rather strict  summability condition as in \eqref{summ cond.}. 
\end{remark}

\noindent At last, we give an example for an upper bound for $\tilde{\tau}_{\boldsymbol{X}}(h,m),$ with $h=0$ for simplicity.

\begin{example}\label{Ex:Cart prod} Let $\boldsymbol{X}=(X_k)$ be the fAR$(1)$ process in Example \ref{Ex: upper bound for L^p errors fAR(1)}, with $\xi=\|\psi\|_{\mathcal{L}} < 1$ and  $\nu_4(\varepsilon_0) < \infty.$ Then, due to the definition of $\tilde{\tau}_{\boldsymbol{X}}(0,m),$ and since $\nu_4(X_k - X^{(k)}_k) \leq \frac{2\nu_4(\varepsilon_0)}{1-\xi}\xi^k$ and $\nu_4(X_1) \leq \frac{\nu_4(\varepsilon_0)}{1-\xi}$ for each $k,$ see the conversions in the proof of Proposition \ref{Prop. L^p-m for AR}, it holds 
\begin{align*}
\tilde{\tau}_{\boldsymbol{X}}(0,m) &\leq \frac{\nu^4_4(\varepsilon_0)}{(1-\xi)^4}\left[\,1 + \frac{8}{1-\xi}\cdot\begin{cases}
 \,\xi, &\textnormal{if}\;m\in\{1,2\},\\
 \,\xi^{m-1}/(2m-3), &\textnormal{if}\;m>2,
  \end{cases}\right].
\end{align*}
%Thereby, for instance for $h=1, \xi=\frac{1}{2},$ we have $\tau_{\boldsymbol{X}}(1) \leq 80\,\nu^4_4(\varepsilon_0)$
Consequently for $m=m_{N}\rightarrow\infty$ with $m=\mathrm o(N)$ holds
\begin{align*}
\limsup_{N\rightarrow\infty}\,\tilde{\tau}_{\boldsymbol{X}}(0,m) &\leq \frac{\nu^4_4(\varepsilon_0)}{(1-\xi)^4}\,.
\end{align*}
\end{example}

\vspace{-.05in}

\section{Consequences of our results}\label{Section: Consequences}

Our upper bounds allow to determine a minimum sample size $N$ so that a certain threshold $\tau > 0$ is not exceeded by the estimation errors,  e.g.,  $\Erw\!\|\hat{\mathscr{C}}_{\!\boldsymbol{X}^{[m]}}\! - \mathscr{C}_{\!\boldsymbol{X}^{[m]}}\|^2_{\mathcal{S}}\leq \tau$ if $\frac{m^2}{N-m+1}\,\tilde{\tau}_{\boldsymbol{X}}(0,m) \leq \tau.$ Although the values $\xi_{\boldsymbol{X}}(h,m,n),$ $\tau_{\boldsymbol{X}}(h,m,n)$ and thus $\tilde{\tau}_{\boldsymbol{X}}(h,m)$ are generally not known, in certain situations they are bounded above by expressions that are known or can be chosen (e.g. \frenchspacing when simulating), see Example \ref{Ex:Cart prod}. Further, our results can be used to determine upper bounds of estimation errors for covariance operator's eigenelements, parameters of the well-known fAR(MA) models, and certain sums and products involving lagged (cross-)
covariance operators, which will be discussed below. 

\vspace{-.025in}

\subsection{Consequences for eigenelements}\label{Subsection: Consequences for eigenelements}

Let $(\hat{\lambda}_j, \hat{c}_j)_{j\in\mathbb{N}}=(\hat{\lambda}_j(m), \hat{c}_j(m))_{j\in\mathbb{N}}$ and $(\lambda_j, c_j)_{j\in\mathbb{N}}=(\lambda_j(m), c_j(m))_{j\in\mathbb{N}}$ be the sequences of eigenelement pairs, i.e. \frenchspacing pairs of eigenvalues associated with the eigenfunctions, of the covariance operators $\hat{\mathscr{C}}_{\!\boldsymbol{X}^{[m]}}$ and $\mathscr{C}_{\!\boldsymbol{X}^{[m]}},$ respectively. Corollary \ref{Cor cov op} has several consequences for the eigenelements discussed here, where we always assume $m=m_N\rightarrow\infty$ with $m=\mathrm o(N).$

Due to $\sup_{j\in\mathbb{N}}|\hat{\lambda}_j - \lambda_j| \leq \|\hat{\mathscr{C}}_{\!\boldsymbol{X}^{[m]}}\! -  \mathscr{C}_{\!\boldsymbol{X}^{[m]}}\|_{\mathcal{L}}$ for any $m$ after \cite{Bosq2000}, Lemma 4.2, and  $\|\cdot\|_{\mathcal{L}} \leq \|\cdot\|_{\mathcal{S}},$ it holds 
\begin{align*}
	\limsup_{N\rightarrow\infty}\frac{N}{m^3}\Erw\sup_{j\in\mathbb{N}}|\hat{\lambda}_j - \lambda_j|^2 \leq \nu^4_4(X_1).   
\end{align*}
If the eigenspace associated to the eigenvalues $\lambda_j$ is one-dimensional, where $\lambda_1 \geq \lambda_2 \geq \cdots \geq 0$ can be assumed without loss of generality (w.l.o.g.), the estimates $\hat{c}'_j \coloneqq \sgn\langle\hat{c}_j, c_j\rangle \hat{c}_j$ for $c_j,$ with  $\sgn(x)\coloneqq \mathds 1_{[0, \infty)}(x) - \mathds 1_{(-\infty, 0)}(x), \,x \in \mathbb{R},$ satisfy  $\|\hat{c}'_j - c_j\| \leq \frac{2\sqrt{2}}{\alpha_j}\|\hat{\mathscr{C}}_{\!\boldsymbol{X}^{[m]}}\! -  \mathscr{C}_{\!\boldsymbol{X}^{[m]}}\|_{\mathcal{L}}$ for each $j,m$ after \cite{Bosq2000}, Lemma 4.3, where  $\alpha_1=\alpha_1(m)\coloneqq \lambda_1 - \lambda_2,$ and $\alpha_j=\alpha_j(m)\coloneqq \min\{\lambda_{j-1} - \lambda_j,\lambda_j - \lambda_{j+1}\}$ for $j>1.$ Then, provided  $\limsup_{N\rightarrow\infty}\alpha^{-1}_j(m) \leq a_j$ for some $a_j\in [0,\infty)$ being  independent of $m,$ 
\begin{gather}
	\limsup_{N\rightarrow\infty}\frac{N}{m^3}\Erw\!\|\hat{c}'_j - c_j\|^2 \leq 8\,a^2_j\,\nu^4_4(X_1),\quad j\in\mathbb{N}. \label{classical upper bound for eigenfunctions}
\end{gather}
\cite{Bosq2000}, Lemma 4.3 also gives $\sup_{1\leq j \leq k}\|\hat{c}'_j - c_j\| \leq 2\sqrt{2}\Lambda_k\|\hat{\mathscr{C}}_{\!\boldsymbol{X}^{[m]}}\! -  \mathscr{C}_{\!\boldsymbol{X}^{[m]}}\|_{\mathcal{L}}$ for each $k\in\mathbb{N},$ where  $\Lambda_k = \Lambda_k(m) \coloneqq \sup_{1\leq j \leq k}(\lambda_j - \lambda_{j+1})^{-1}.$ Thus, for each $k,$ provided $\limsup_{N\rightarrow\infty}\Lambda_k(m) \leq L_k$ for some $L_k\in[0,\infty)$ which is independent of $m,$ it holds 
\begin{gather}
	\limsup_{N\rightarrow\infty}\frac{N}{m^3}\Erw\!\sup_{1\leq j \leq k}\!\|\hat{c}'_j - c_j\|^2 \leq 8\,L^2_k\,\nu^4_4(X_1),\quad k\in\mathbb{N}. 
\end{gather}
The variable $k$ can also depend on $N$ which is needed when estimating complete parameters in fTS models by projecting onto a $k$-dimensional subspace and letting $k=k_N$ tend to infinity, see \cite{Bosq2000}, \cite{Kuehnert2022}, \cite{KuehnertRiceAue2024+}, \cite{kuenzer:2024}, to name a few. For $k\rightarrow\infty,$  we get
\begin{align}\label{uniform eigenfunction estimation}
	\Erw\!\sup_{1\leq j \leq k}\!\|\hat{c}'_j - c_j\|^2 = \mathrm O(\Lambda^{2}_k\,m^3N^{-1}).
\end{align}

In the following, we analyze the relationship between the eigenvalues of the covariance operators of the Cartesian product space-valued to the original process in certain situations. Such a relationship cannot be derived without assuming a certain structure, but we know that
\begin{align}
	\sum^\infty_{j=1}\,\lambda_j(m) = \|\mathscr{C}_{\!\boldsymbol{X}^{[m]}} \|_{\mathcal{N}} %\Erw\!\|X^{[m]}_0\|^2 = m\Erw\!\|X_0\|^2 
	= m\|\mathscr{C}_{\!\boldsymbol{X}}\|_{\mathcal{N}}= m\sum^\infty_{j=1}\,\lambda_j(1), \quad m\in\mathbb{N},
\end{align}
where $\lambda_j(1)$ denote the eigenvalues of $\mathscr{C}_{\!\boldsymbol{X}^{[1]}}=\mathscr{C}_{\!\boldsymbol{X}}.$

\begin{example} Throughout, let $\boldsymbol{X}=(X_k)_{k\in\mathbb{Z}}\subset\mathcal{H}$ be a process which   satisfies Assumption \ref{As: X L^2, centered}. 
	\begin{itemize}
		\item[\textnormal{(a)}] Suppose $(X_k)$ is i.i.d. Then $\mathscr{C}_{\!\boldsymbol{X}^{[m]}}=\mathrm{diag}(\mathscr{C}_{\!\boldsymbol{X}}, \dots, \mathscr{C}_{\!\boldsymbol{X}})$ for each $m,$ and the eigenvalues satisfy 
		$$\underbrace{\lambda_1(m) = \cdots = \lambda_m(m)}_{=\,\lambda_1(1)} \,\geq\, \underbrace{\lambda_{m+1}(m) = \cdots = \lambda_{2m}(m)}_{=\, \lambda_2(1)} \,\geq\, \underbrace{\lambda_{2m+1}(m) = \cdots = \lambda_{3m}(m)}_{=\,\lambda_3(1)} \,\geq\, \cdots. $$ 
		Hence, \eqref{classical upper bound for eigenfunctions}--\eqref{uniform eigenfunction estimation} cannot  be used directly for $m>1.$ For upper bounds considering multi-dimensional eigenspaces (with fixed dimension), we refer to \cite{Bosq2000}, Lemma 4.4.
		\item[\textnormal{(b)}] Let $(X_k)\subset L^2_\mathcal{H}$ be  the fAR process in \eqref{equation fAR(1)} with self-adjoint fAR operator satisfying $\psi(c_j(1)) = \psi_jc_j(1)$ for each $j$ where  $1>\psi_1>\psi_2 > \cdots > 0.$ For the entries of $\mathscr{C}_{\!\boldsymbol{X}^{[m]}}=(\mathscr{C}^{j-i}_{\!\boldsymbol{X}})_{1\leq i, j \leq m}$ holds  $\mathscr{C}^{j-i}_{\!\boldsymbol{X}} = \psi^{j-i}\mathscr{C}_{\!\boldsymbol{X}}.$  % thus, for each $m$ holds 
		%$\mathscr{C}_{\!\boldsymbol{X}^{[m]}} = (\psi^{|j-i|}\mathscr{C}_{\!\boldsymbol{X}})_{1\leq i,j\leq m},$ The eigenvalues of $\mathscr{C}_{\!\boldsymbol{X}}=\mathscr{C}_{\!\boldsymbol{X}^{[1]}},$ $\mathscr{C}_{\!\boldsymbol{X}^{[2]}},$ $ \mathscr{C}_{\!\boldsymbol{X}^{[3]}},$ $\dots, $$ \mathscr{C}_{\!\boldsymbol{X}^{[m]}}, m\in\mnathbb{N},$ are $\lambda_k(1),$ $\lambda_k(2) = (1+\psi_k)\lambda_k(1),$ $\lambda_k(3) = (....)\lambda_k(1), \dots, \lambda_k(m) = (....)\lambda_k(1),$ respectively. In fact, for each given eigenfunction $k,$ there are even more eigenvalues, e.g., $(1\pm\psi_k)\lambda_k(1)$ for $m=2,$ but since we order the  eigenvalues in a decreasing manner w.l.o.g., only the largest eigenvalues are relevant for us. 
		We also assume injectivity of $\mathscr{C}_{\boldsymbol{\varepsilon}}$, which transfers to $\mathscr{C}_{\! \boldsymbol{X}^{[m]}}$ for each $m,$ see \cite{KuehnertRiceAue2024+}, Proposition 3.1, so  all eigenvalues $\lambda_j(m)$ of $\mathscr{C}_{\!\boldsymbol{X}^{[m]}}$ are strictly positive. The explicit determination of all $\lambda_j(m),$ and their relation to the eigenvalues $\lambda_j(1)$ associated with the eigenfunctions $c_j(1)$ of $\mathscr{C}_{\! \boldsymbol{X}}$ is difficult, maybe even impossible, but with the normed vector $\tilde{c}_j(m)\coloneqq\frac{1}{\sqrt{m}}(c_j(1), \dots, c_j(1))^\top\in\mathcal{H}^m,$ we can derive an upper bound for the largest of the $m$ eigenvalues for each $j.$ We have
		\begin{align*}
			\mathscr{C}_{\!\boldsymbol{X}}(\tilde{c}_j(m))&=\frac{\lambda_j}{\sqrt{m}}\Big(\sum^m_{j=1}\psi^{|j-1|}_jc_j(1), \sum^m_{j=1}\psi^{|j-2|}_jc_j(1), \dots, \sum^m_{j=1}\psi^{|j-m|}_jc_j(1)\Big)^{\!\top}\\
			&= \frac{\lambda_j}{\sqrt{m}(1-\psi_j)}\Big(\psi_{j,1,m}c_j(1), \psi_{j,2,m}c_j(1), \dots, \psi_{j,m,m}c_j(1)\Big)^{\!\top},
		\end{align*}
		with $\psi_{j,\ell,m}\coloneqq 1+\psi_j(1-(\psi^{\ell-1}_j+\psi^{m-\ell}))$ for each $j,m$ and $\ell=1, \dots, m,$  where   
		\begin{align*}
			\max_{1\leq\ell\leq m}\psi_{j,\ell,m}&=\psi_{j,\lfloor\frac{m+1}{2}\rfloor, m} = 1+\psi_j\big(1-2\psi^{\lfloor\frac{m+1}{2} \rfloor-1}_k\big)\,.
		\end{align*}
		Thus, for each $m,$ the eigenvalues $\lambda_j(m)$ satisfy 
		\begin{align*}
			\lambda_j(m) \leq \frac{1+\psi_j\big(1-2\psi^{\lfloor\frac{m+1}{2} \rfloor-1}_j\big)}{1-\psi_j}\,\lambda_j(1)\,\stackrel{m\rightarrow\infty}{\longrightarrow}\,\frac{1+\psi_j}{1-\psi_j}\,\lambda_j(1),
		\end{align*}
		and if the operator $\psi$ is compact, $\lambda_j(m)\rightarrow 0$ for $j\rightarrow\infty$ both for fixed $m$ and for $m\rightarrow\infty.$
		%To do further analysis, though, let  $\psi\coloneqq \frac{1}{2}(c_1(1)\otimes c_2(1)),$ where $c_j(1)$ are the eigenfunctions of $\mathscr{C}_{\!\boldsymbol{X}}$ associated to the eigenvalues $\lambda_j(1).$ Then, $\psi^\ast= \frac{1}{2}(c_2(1)\otimes c_1(1)),$ where $\psi^i=\psi^{\ast i}$ for $i>1,$ thus the Toeplitz matrix $\mathscr{C}_{\!\boldsymbol{X}^{[m]}}$ is tridiagonal and has the form 
		%  \begin{align}\label{Operator-valued covop-matrix in fAR(1)}
			%    \mathscr{C}_{\!\boldsymbol{X}^{[m]}} = \begin{pmatrix}
				%      \mathscr{C}_{\!\boldsymbol{X}}              & \psi\mathscr{C}_{\!\boldsymbol{X}}                             & 0 & \cdots & 0\\
				%      \mathscr{C}_{\!\boldsymbol{X}}\psi^{{\ast}} & \mathscr{C}_{\!\boldsymbol{X}} & \psi\mathscr{C}_{\!\boldsymbol{X}} & \ddots & \vdots\\
				%      0 & \mathscr{C}_{\!\boldsymbol{X}}\psi^{{\ast}} & \ddots & \ddots & 0\\
				%            \vdots & \ddots & \ddots & \ddots & \psi\mathscr{C}_{\!\boldsymbol{X}}\\
				%      0 & \cdots & 0 & \mathscr{C}_{\!\boldsymbol{X}}\psi^{{\ast}} & \mathscr{C}_{\!\boldsymbol{X}}
				%    \end{pmatrix}, \qquad m\in\mathbb{N}.
			%  \end{align}  
		%  Due to $\mathscr{C}_{\!\boldsymbol{X}}=\sum_j\lambda_j(1)c_1(1)\!\otimes\!c_1(1),$ it holds $\psi\mathscr{C}_{\!\boldsymbol{X}}(x)=\frac{\lambda_1(1)}{2}(c_1(1)\!\otimes\!c_2(1))(x)$ and $\mathscr{C}_{\!\boldsymbol{X}}\psi^\ast(x)= \frac{\lambda_1(1)}{2}(c_2(1)\otimes c_1(1))(x)$ for any $x\in\mathcal{H}.$
		The upper bounds  \eqref{classical upper bound for eigenfunctions}--\eqref{uniform eigenfunction estimation} can be used here where $\Lambda_j=\Lambda_j(m)$ depends only negligibly on large $m,$ provided the eigenspace of $\lambda_j(m)$ is one-dimensional. This, though,  is not necessarily the case since we have only derived upper bounds for the eigenvalues, and since one of the $m$ eigenvalues derived from $\lambda_j(1)$ could occur again for another $j.$
		\item[\textnormal{(c)}] Now, we suppose the extreme case $X_k = X$ a.s. \frenchspacing for each $k$, where $X\in L^2_\mathcal{H}$ is non-degenerate. Then, 
		$\mathscr{C}_{\!\boldsymbol{X}^{[m]}} = (\mathscr{C}_{\!\boldsymbol{X}})_{1\leq i,j\leq m}$ for each $m.$  $\mathscr{C}_{\!\boldsymbol{X}^{[m]}}$ has the eigenvalues $\lambda_j(m)=m\lambda_j(1)$ for each $j,m$ and 0 (with infinite-dimensional eigenspace), since $\mathscr{C}_{\!\boldsymbol{X}^{[m]}}$ is not injective even if $\mathscr{C}_{\! \boldsymbol{X}}$ is. 
		Furthermore, we assume that the eigenvalues of $\mathscr{C}_{\! \boldsymbol{X}}$ satisfy $\lambda_j(1) = j^{-2}$ for each $j$ (these eigenvalues are up to the factor $\pi^{-2}$ those of the covariance operator of the \emph{Brownian bridge} on the domain $[0,1],$ see \cite{JaimezValderrama1987}). Then \eqref{classical upper bound for eigenfunctions} becomes
		\begin{align*}
			\limsup_{N\rightarrow\infty}\frac{N}{m^{3}}\Erw\!\|\hat{c}'_j - c_j\|^2 \leq 8\,\nu^4_4(X_1)\,\frac{j^2(j+1)^2}{1+2j}\,,\quad j\in\mathbb{N}. 
		\end{align*}
		Further, due to $\Lambda_k = \Lambda_k(m) = \frac{k(k+1)^2}{m(2+1/k)}$ for each $k,$ for \eqref{uniform eigenfunction estimation} holds for $k=k_N\rightarrow\infty\colon$
		\begin{align*}
			\Erw\!\sup_{1\leq j \leq k}\!\|\hat{c}'_j - c_j\|^2 = \mathrm O(k^6mN^{-1}).
		\end{align*}
	\end{itemize}
\end{example}

\vspace{-.025in}

\subsection{Consequences for estimating fAR(MA) parameters}

Consider the  fAR$(p)$ process $\boldsymbol{X}=(X_k)_{k\in\mathbb{Z}}\subset\mathcal{H},$ with $p\in\mathbb{N},$ defined through $X_k = \sum^p_{i=1}\psi_i(X_{k-i}) + \varepsilon_k$ a.s. \frenchspacing for each $k,$ where $(\varepsilon_k)\subset L^2_{\mathcal{H}}$ is i.i.d. \frenchspacing and centered, and where all $\psi_i\in\mathcal{S_H}$ are H-S operators. Then, we have the following \emph{Yule-Walker} equation
\begin{align}\label{Identity pseudo inv}
    \mathscr{C}^1_{\!\boldsymbol{X}^{[p]}, \boldsymbol{X}}\! = \Psi_{\!1}\mathscr{C}_{\!\boldsymbol{X}^{[p]}}
\end{align}
from which $\Psi_{\!1}\coloneqq (\psi_1\, \cdots \,\psi_p)$ can be uniquely estimated if the covariance operator $\mathscr{C}_{\!\boldsymbol{X}^{[p]}}$ is injective. However, \eqref{Identity pseudo inv} is an ill-posed problem because the inverse of $\mathscr{C}_{\!\boldsymbol{X}^{[p]}}$ is not bounded. This problem can be solved, e.g., by using a Moore-Penrose generalized inverse as in  \cite{Reimherr2015} and \cite{kuenzer:2024}, or the Tychonoff-regularization $\hat{\mathscr{C}}^\dagger_{\!\!\boldsymbol{X}^{[p]}}\coloneqq (\hat{\mathscr{C}}_{\!\boldsymbol{X}^{[p]}} + \vartheta_{\!N})^{-1}$ as in \cite{KuehnertRiceAue2024+}, where $(\vartheta_{\!N})_{\!N\in\mathbb{N}}\subset(0,\infty)$ is an appropriate null sequence. For example, we can use the estimate $\hat{\Psi}_{\!1}\coloneqq \hat{\mathscr{C}}^1_{\!\boldsymbol{X}^{[p]}, \boldsymbol{X}}\hat{\mathscr{C}}^\dagger_{\!\! \boldsymbol{X}^{[p]}}$ for $\Psi_{\!1}.$ By projecting this estimate onto the linear space spanned by the eigenfunctions $\hat{c}_1, \dots, \hat{c}_K$ of $\hat{\mathscr{C}}_{\! \boldsymbol{X}^{[p]}}$ for some $K\in\mathbb{N},$ and given that $\Psi_{\!1}$ can be approximated sufficiently fast in the sense of a \emph{Sobolev condition} with parameter $\beta>0,$ based on Proposition \ref{Corollary Cross-cov on X's}, for an appropriately chosen sequence $K=K_N\rightarrow\infty$ holds
\begin{align*}
  \|\hat{\Psi}_{\!1} - \Psi_{\!1}\|_{\mathcal{S}} = \mathrm O_{\Prob}(K^{-\beta}),
\end{align*}
 see Section 4.1 in \cite{KuehnertRiceAue2024+}. Note that although we have upper bounds for the estimation errors of the lagged covariance and cross-covariance operators in the $L^2$-sense, we can only give upper bounds of $\|\hat{\Psi}_{\! 1} - \Psi_{\!1}\|_{\mathcal{S}}$ in the sense of stochastic boundedness, since these upper bounds contain reciprocals of eigenvalue estimates which only converge stochastically.

The following example illustrates why allowing the Cartesian power(s) to go to infinity is of great benefit in parameter estimation of certain fTS models. Consider a fARMA$(p,q)$ process $\boldsymbol{X}=(X_k)_{k\in\mathbb{Z}}\subset\mathcal{H},$ with $p,q\in\mathbb{N},$ i.e \frenchspacing $X_k = \sum^p_{i=1}\alpha_i(X_{k-i}) + \sum^q_{j=1}\beta_j(\varepsilon_{k-j}) + \varepsilon_k$ a.s. \frenchspacing for any $k,$ where $(\varepsilon_k)_{k\in\mathbb{Z}}\subset L^2_{\mathcal{H}}$ is i.i.d. \frenchspacing and centered, and $\alpha_i, \beta_j\in\mathcal{S_H}$ are assumed to be H-S operators. If such a process is a fAR$(\infty)$ process (or equivalently \emph{invertible}) with finite second moments, it satisfies the 'approximate' Yule-Walker equation
\begin{align}\label{Identity pseudo inv fARMA}
    \mathscr{C}^1_{\!\boldsymbol{X}^{[m]}, \boldsymbol{X}}\! = \Psi_{\!2}\mathscr{C}_{\!\boldsymbol{X}^{[m]}} + \sum_{\ell>m}\psi_\ell\mathscr{C}^{1-\ell}_{\!\boldsymbol{X}^{[m]}, \boldsymbol{X}},\quad m\in\mathbb{N},
\end{align}
where $\Psi_{\!2} = \Psi_{\!2}(m)\coloneqq (\psi_1\, \cdots \,\psi_m).$ Provided that all $\psi_i$ are H-S operators with $\sum_{i\in\mathbb{N}}\|\psi_i\|^2_{\mathcal{S}} < \infty$ the series in \eqref{Identity pseudo inv fARMA} exists for every $m,$ and vanishes for $m\rightarrow\infty.$ Moreover, under some 'uniform' Sobolev condition with $\beta>0$ and other technical conditions, for $\hat{\Psi}_{\!2}\coloneqq \hat{\mathscr{C}}^1_{\! \boldsymbol{X}^{[m]}, \boldsymbol{X}}\hat{\mathscr{C}}^\dagger_{\!\! \boldsymbol{X}^{[m]}}$ projected onto a certain $K$-dimensional space (as above) holds with appropriate sequences $m=m_N, K=K_N\rightarrow\infty$ and by using Proposition \ref{Corollary Cross-cov on X's}:
\begin{align}\label{result for invertible representation}
  \|\hat{\Psi}_{\!2} - \Psi_{\!2}\|_{\mathcal{S}} = \mathrm O_{\Prob}(K^{-\beta}).
\end{align}
\noindent In turn, if $(X_k)$ is causal and linear,  fARMA$(p,q)$ can be expressed by fAR$(\infty)$ processes, thus based on \eqref{result for invertible representation} and similar technical conditions as above, consistency results can also be derived for the fARMA$(p,q)$ operators, see \cite{KuehnertRiceAue2024+} for a detailed explanation.

\begin{remark} Our upper bounds are also useful in the context of estimating f(G)ARCH parameters, since the estimation procedure is closely related to estimating fAR(MA) parameters, see \cite{Hoermannetal2013} and \cite{Kuehnert2020}. 
\end{remark}

\vspace{-.025in}

\subsection{Consequences for parameter estimation in sums and products}

Here we illustrate more abstractly the consequences of our results for parameter estimation in certain situations. Let $\Delta$ be a parameter that can be consistently estimated by an estimator $\hat{\Delta}$ based on the samples $X_1, \dots, X_N$ and possibly also $Y_1, \dots, Y_N$ of the processes $\boldsymbol{X}=(X_k)\subset\mathcal{H}$ and $\boldsymbol{Y}=(Y_k)\subset\mathcal{H}_\star$ in Section \ref{Section: Main results}, respectively. We also assume that there is a function $f\colon\mathbb{N}\times \Theta\rightarrow [0,\infty)$ with $f(N,\boldsymbol{\theta})\rightarrow\infty,$ where $\boldsymbol{\theta}\coloneqq (\theta_1, \dots, \theta_d)^\top\in\Theta\subset \mathbb{R}^d,$ with $d\in\mathbb{N},$ are parameters which may depend on $N,$ such that for some $\tau\geq 0$ holds
\begin{align}\label{Inequ Appl}
  f(N,\boldsymbol{\theta})\Erw\!\|\hat{\Delta} - \Delta\|^2 \leq \tau,\quad  N\in\mathbb{N},\, \boldsymbol{\theta}\in\Theta,   
\end{align}
 and where the norm $\|\cdot\|$ is given below. For example, in the context of estimating our lagged covariance operators $\mathscr{C}^h_{\!\boldsymbol{X}^{[m]}},$ we have $\boldsymbol{\theta} = (h, m)^\top$ and $f(N,\boldsymbol{\theta})=\frac{N-|h|-m+1}{m^2}.$
 
First, consider that for some $\Psi_{\!3}$ and $\Delta \in\mathcal{S}_{\mathcal{H}^m, \mathcal{H}^n_{\star}}$ holds with $h\in\mathbb{Z},$ $m,n\in\mathbb{N}\colon$
\begin{align*}
    \Psi_{\!3} = \Delta + \mathscr{C}^h_{\!\boldsymbol{X}^{[m]}, \boldsymbol{Y}^{[n]}}.
\end{align*}
By using the plug-in estimator $\hat{\Psi}_{\!3}\coloneqq \hat{\Delta} + \hat{\mathscr{C}}^h_{\!\boldsymbol{X}^{[m]}, \boldsymbol{Y}^{[n]}}$ for $\Psi_{\!3}\in\mathcal{S}_{\mathcal{H}^m, \mathcal{H}^n_{\star}},$ according to Theorem \ref{Theo cross-cov op}, \eqref{Inequ Appl} and elementary conversions, we obtain  for each $h,m,n, \boldsymbol{\theta},$ and with $\kappa'$ in Theorem \ref{Theo cross-cov op}:
\begin{align}
  \Erw\!\|\hat{\Psi}_{\!3} - \Psi_{\!3}\|^2_{\mathcal{S}} %&\leq 2 \Erw\!\|\hat{\Delta} - \Delta\|^2_{\mathcal{S}} + 2 \Erw\!\|\hat{\mathscr{C}}^h_{\!\boldsymbol{X}^{[m]}, \boldsymbol{Y}^{[n]}} - \mathscr{C}^h_{\!\boldsymbol{X}^{[m]}, \boldsymbol{Y}^{[n]}}\|^2_{\mathcal{S}}\notag\\
  &\leq \!\frac{2\tau}{f(N,\boldsymbol{\theta})} +  \frac{2mn(2\kappa'-1)}{N'}\,\xi_{\boldsymbol{X,Y}}(h, m, n).\notag
\end{align}

Now suppose for some  $\Psi_{\!4}$ and $\Delta \in\mathcal{T}_{\mathcal{H}^n_{\star}, \mathcal{H}^p_{\star\star}},$ with $\mathcal{T}\in\{\mathcal{L, S, N}\}$ and $n,p\in\mathbb{N},$ and where $\mathcal{H}_{\star\star}$ is another real, separable Hilbert space, holds for certain $h\in\mathbb{Z},$ $m,n\in\mathbb{N}\colon$ 
\begin{align*}
    \Psi_{\!4} = \Delta\mathscr{C}^h_{\!\boldsymbol{X}^{[m]}, \boldsymbol{Y}^{[n]}}\,.
\end{align*}
Then, $\hat{\Psi}_{\!4}\coloneqq \hat{\Delta} \hat{\mathscr{C}}^h_{\!\boldsymbol{X}^{[m]}, \boldsymbol{Y}^{[n]}}$ is a reasonable estimator for $\Psi_{\!4}\in\mathcal{S}_{\mathcal{H}^m\!, \mathcal{H}^p_{\star\star}},$ and due to triangle, operator-valued H\"older's, Jensen's, and Cauchy-Schwarz inequality (CSI), Theorem \ref{Theo cross-cov op},  \eqref{Inequ Appl}, and  also $\|\mathscr{C}^h_{\!\boldsymbol{X}^{[m]}, \boldsymbol{Y}^{[n]}}\! \|_{\mathcal{S}}\leq \sqrt{mn}\,\nu_2(X_0)\nu_2(Y_0),$ we get for each $h,m,n, \boldsymbol{\theta}\colon$  
\begin{align*}
\Erw\!\|\hat{\Psi}_{\!4} - \Psi_{\!4}\|_{\mathcal{S}} %&\leq \Erw\!\|\hat{\Delta} - \Delta\|_{\mathcal{L}}\big\|\hat{\mathscr{C}}^h_{\!\boldsymbol{X}^{[m]}, \boldsymbol{Y}^{[n]}}\! \big\|_{\mathcal{S}} + \|\Delta\|_{\mathcal{L}}\Erw\!\big\|\hat{\mathscr{C}}^h_{\!\boldsymbol{X}^{[m]}, \boldsymbol{Y}^{[n]}}\! - \mathscr{C}^h_{\!\boldsymbol{X}^{[m]}, \boldsymbol{Y}^{[n]}}\!\big\|_{\mathcal{S}}\notag\\
&\leq \Erw\!\Big[\,\big\|\hat{\mathscr{C}}^h_{\!\boldsymbol{X}^{[m]}, \boldsymbol{Y}^{[n]}}\! - \mathscr{C}^h_{\!\boldsymbol{X}^{[m]}, \boldsymbol{Y}^{[n]}}\!\big\|_{\mathcal{S}}\big(\|\hat{\Delta} - \Delta\|_{\mathcal{L}} + \|\Delta\|_{\mathcal{L}} \big)\Big] + \big\|\mathscr{C}^h_{\!\boldsymbol{X}^{[m]}, \boldsymbol{Y}^{[n]}}\! \big\|_{\mathcal{S}}\Erw\!\|\hat{\Delta} - \Delta\|_{\mathcal{L}}\notag\\
&\leq \sqrt{\frac{mn(2\kappa'-1)\,\xi_{\boldsymbol{X,Y}}(h, m, n)}{N'}}\,\bigg(\sqrt{\frac{\tau}{f(N,\boldsymbol{\theta})}}\, + \|\Delta\|_{\mathcal{L}}\bigg) + \sqrt{\frac{mn\tau}{f(N,\boldsymbol{\theta})}}\,\nu_2(X_0)\nu_2(Y_0).
\end{align*}
Moreover, if \eqref{Inequ Appl} holds with the H-S norm $\|\cdot\|_{\mathcal{S}},$ the latter upper bound with  $\|\Delta\|_{\mathcal{L}}$ replaced by $\|\Delta\|_{\mathcal{S}}$ is an upper bound for  $\Erw\!\|\hat{\Psi}_{\!4} - \Psi_{\!4}\|_{\mathcal{N}}.$ 

\vspace{-.05in}

\section{Concluding remarks}\label{Conclusion}

This article derives upper bounds of the estimation errors for lagged cross-covariance operators of Cartesian product Hilbert space-valued processes defined by two $L^p$-$m$-approximable processes in potentially different general separable Hilbert spaces. The upper bounds are stated in the $L^2$-sense for any sample size $N,$ lag $h$ and Cartesian powers $m$ and $n,$ which makes them novel in all their generality. A slightly sharper bound is obtained when the two processes that define the Cartesian product space-valued processes are (a.s.) identical. We also show how our upper bounds can be used to estimate eigenelements. The upper bounds for the eigenvalues coincide with those for the covariance operators, but the situation for the eigenfunctions is somewhat more delicate due to the appearance of reciprocal spectral gaps, and because both the decay behavior of the eigenvalues and the relation to the Cartesian power tending to infinity may need to be controlled. The examples show that a strong dependence structure has a positive impact on the upper bounds. Furthermore, we show that our upper bounds for lagged (cross-)covariance operators of Cartesian product Hilbert space-valued processes are used in combination with other technical conditions to consistently estimate fAR$(m)$ operators for fixed $m$ as well as for fARMA$(p,q)$ operators, rewritten by fAR$(\infty)$, where letting $m$ tend to infinity is required. 

We would also like to point out that possible further research topics include the study of the asymptotic distribution of the estimation errors, and the application of all our results to general Banach or even metric spaces. 

\vspace{-.05in}

\section*{Acknowledgements}

The author thanks the co-editor and two referees for their careful reading, questions, and thoughtful insights, which greatly improved the quality of the paper. Thanks also go to Alexander Aue, Xiucai Ding, Miles Lopes (University of California, Davis), Gregory Rice (University of Waterloo), and Colin Decker (Ruhr-Universität Bochum) for valuable general discussions and comments. The author would also like to point out that most of the work on this article was done while he was at University of California, Davis, and the final version was written at Ruhr-Universität Bochum.  

\newpage

\bibliography{Estimation_Lagged_Lpm}

\begin{thebibliography}{}

\bibitem[\protect\citeauthoryear{Aue, Horv\'ath, and Pellatt}{Aue
  et~al.}{2017}]{Aueetal2017}
Aue, A., L.~Horv\'ath, and D.~Pellatt (2017).
\newblock Functional generalized autoregressive conditional heteroskedasticity.
\newblock {\em Journal of Time Series Analysis\/}~{\em 38}, 3--21.

\bibitem[\protect\citeauthoryear{Aue and Klepsch}{Aue and
  Klepsch}{2017}]{AueKlepsch2017}
Aue, A. and J.~Klepsch (2017).
\newblock Estimating functional time series by moving average model fitting.
  ar{X}iv:1701.00770.
\newblock \url{https://arxiv.org/abs/1701.00770}.

\bibitem[\protect\citeauthoryear{Bosq}{Bosq}{2000}]{Bosq2000}
Bosq, D. (2000).
\newblock {\em Linear processes in function spaces}.
\newblock Lecture notes in statistics. New York: Springer.

\bibitem[\protect\citeauthoryear{Brockwell and Davis}{Brockwell and
  Davis}{1991}]{BrockwellDavis1991}
Brockwell, P. and R.~Davis (1991).
\newblock {\em Time series: theory and methods\/} (2 ed.).
\newblock New York: Springer.

\bibitem[\protect\citeauthoryear{Dedecker, Doukhan, Lang, Jos{\'e}~Rafael,
  Louhichi, and Prieur}{Dedecker et~al.}{2007}]{Dedecker2007}
Dedecker, J., P.~Doukhan, G.~Lang, L.~R. Jos{\'e}~Rafael, S.~Louhichi, and
  C.~Prieur (2007).
\newblock {\em Weak dependence: with examples and applications}.
\newblock Lecture Notes in Statistics. New York: Springer.

\bibitem[\protect\citeauthoryear{Ferraty and Vieu}{Ferraty and
  Vieu}{2006}]{FerratyVieu2006}
Ferraty, F. and P.~Vieu (2006).
\newblock {\em Nonparametric functional data analysis}.
\newblock New York: Springer.

\bibitem[\protect\citeauthoryear{H\"ormann, Horv\'ath, and Reeder}{H\"ormann
  et~al.}{2013}]{Hoermannetal2013}
H\"ormann, S., L.~Horv\'ath, and R.~Reeder (2013).
\newblock A functional version of the {ARCH} model.
\newblock {\em Econometric Theory\/}~{\em 29\/}(2), 267--288.

\bibitem[\protect\citeauthoryear{H\"ormann and Kidzi\'nski}{H\"ormann and
  Kidzi\'nski}{2015}]{HoermannKidzinski2015}
H\"ormann, S. and L.~Kidzi\'nski (2015).
\newblock A note on estimation in hilbertian linear models.
\newblock {\em Scandinavian Journal of Statistics\/}~{\em 42\/}(1), 43--62.

\bibitem[\protect\citeauthoryear{H\"ormann and Kokoszka}{H\"ormann and
  Kokoszka}{2010}]{HoermannKokoszka2010}
H\"ormann, S. and P.~Kokoszka (2010).
\newblock Weakly dependent functional data.
\newblock {\em Annals of Statistics\/}~{\em 38}, 1845--1884.

\bibitem[\protect\citeauthoryear{Horv\'ath and Kokoszka}{Horv\'ath and
  Kokoszka}{2012}]{HorvathKokoszka2012}
Horv\'ath, L. and P.~Kokoszka (2012).
\newblock {\em Inference for functional data with applications}.
\newblock New York: Springer.

\bibitem[\protect\citeauthoryear{Hsing and Eubank}{Hsing and
  Eubank}{2015}]{HsingEubank2015}
Hsing, T. and R.~Eubank (2015).
\newblock {\em Theoretical foundations of functional data analysis, with an
  introduction to linear operators}.
\newblock West Sussex: Wiley.

\bibitem[\protect\citeauthoryear{Jaimez and Valderrama}{Jaimez and
  Valderrama}{1987}]{JaimezValderrama1987}
Jaimez, R. and M.~Valderrama (1987).
\newblock On the {K}arhunen-{L}oève expansion for transformed processes.
\newblock {\em Trabajos de Estadistica\/}~{\em 2\/}(2), 81--90.

\bibitem[\protect\citeauthoryear{Kokoszka and Reimherr}{Kokoszka and
  Reimherr}{2017}]{kokoszka:2017:FDA-book}
Kokoszka, P. and M.~Reimherr (2017).
\newblock {\em Introduction to functional data analysis}.
\newblock New York: Chapman and Hall/CRC.

\bibitem[\protect\citeauthoryear{Kuenzer}{Kuenzer}{2024}]{kuenzer:2024}
Kuenzer, T. (2024).
\newblock {Estimation of functional {ARMA} models}.
\newblock {\em Bernoulli\/}~{\em 30\/}(1), 117--142.

\bibitem[\protect\citeauthoryear{K\"uhnert}{K\"uhnert}{2020}]{Kuehnert2020}
K\"uhnert, S. (2020).
\newblock Functional {ARCH} and {GARCH} models: {A} {Y}ule-{W}alker approach.
\newblock {\em Electronic Journal of Statistics\/}~{\em 14\/}(2), 4321--4360.

\bibitem[\protect\citeauthoryear{K\"uhnert}{K\"uhnert}{2022}]{Kuehnert2022}
K\"uhnert, S. (2022).
\newblock Lagged covariance and cross-covariance operators of processes in
  {C}artesian products of abstract {H}ilbert spaces.
\newblock {\em Electronic Journal of Statistics\/}~{\em 16\/}(2), 4823--4862.

\bibitem[\protect\citeauthoryear{K\"uhnert, Rice, and Aue}{K\"uhnert
  et~al.}{2024}]{KuehnertRiceAue2024+}
K\"uhnert, S., G.~Rice, and A.~Aue (2024+).
\newblock Estimating invertible processes in {H}ilbert spaces, with
  applications to functional {ARMA} processes. \texttt{ar{X}iv:2407.12221}.
\newblock \url{https://https://arxiv.org/pdf/2407.12221}.

\bibitem[\protect\citeauthoryear{Kutta}{Kutta}{2024}]{Kutta2024+}
Kutta, T. (2024+).
\newblock Approximately mixing time series.
\newblock \url{https://papers.ssrn.com/sol3/papers.cfm?abstract_id=4882128}.

\bibitem[\protect\citeauthoryear{Liebl}{Liebl}{2013}]{Liebl2013}
Liebl, D. (2013).
\newblock {Modeling and forecasting electricity spot prices: {A} functional
  data perspective}.
\newblock {\em Annals of Applied Statistics\/}~{\em 7\/}(3), 1562--1592.

\bibitem[\protect\citeauthoryear{Panaretos and Tavakoli}{Panaretos and
  Tavakoli}{2013}]{panaretos:tavakoli:2013}
Panaretos, V.~M. and S.~Tavakoli (2013).
\newblock Fourier analysis of stationary time series in function space.
\newblock {\em The Annals of Statistics\/}~{\em 41\/}(2), 568--603.

\bibitem[\protect\citeauthoryear{Ramsay and Silverman}{Ramsay and
  Silverman}{2005}]{RamsaySilverman2005}
Ramsay, J. and B.~Silverman (2005).
\newblock {\em Functional data analysis}.
\newblock Springer Series in Statistics. New York: Springer.

\bibitem[\protect\citeauthoryear{Reimherr}{Reimherr}{2015}]{Reimherr2015}
Reimherr, M. (2015).
\newblock Functional regression with repeated eigenvalues.
\newblock {\em Statistics \& Probability Letters\/}~{\em 107}, 62--70.

\bibitem[\protect\citeauthoryear{Rice and Shum}{Rice and
  Shum}{2019}]{RiceShum2019}
Rice, G. and M.~Shum (2019).
\newblock Inference for the lagged cross-covariance operator between functional
  time series.
\newblock {\em Journal of Time Series Analysis\/}~{\em 40\/}(5), 665--692.

\bibitem[\protect\citeauthoryear{Rice, Wirjanto, and Zhao}{Rice
  et~al.}{2020}]{Riceetal2020}
Rice, G., T.~Wirjanto, and Y.~Zhao (2020).
\newblock Forecasting value at risk with intra-day return curves.
\newblock {\em International journal of forecasting\/}~{\em 36\/}(3),
  1023--1038.

\bibitem[\protect\citeauthoryear{Shumway and Stoffer}{Shumway and
  Stoffer}{2017}]{ShumwayStoffer2017}
Shumway, R. and D.~Stoffer (2017).
\newblock {\em Time series analysis and its applications with {R} examples\/}
  (4 ed.).
\newblock New York: Springer.

\bibitem[\protect\citeauthoryear{Stout}{Stout}{1974}]{Stout1974}
Stout, W.~F. (1974).
\newblock {\em Almost sure convergence\/} (2 ed.).
\newblock Springer Series in Statistics. New York: Academic Press.

\end{thebibliography}

\appendix

\vspace{-.05in}

\section{Proofs}\label{Section Proofs}

\begin{proof}[\textbf{Proof of Proposition \ref{Prop. L^p-m for AR}}] The proof uses arguments from \cite{HoermannKokoszka2010}, Example 2.1. Note that our process satisfies $X_k=\sum^\infty_{j=0}\psi^j(\varepsilon_{k-j})$ a.s. \frenchspacing for each $k,$ with $\psi^0\coloneqq \mathbb{I}. $ Furthermore, since $\|\psi^{j_0}\|_{\mathcal{L}} < 1$ for some $j_0\in\mathbb{N}$ is equivalent to that there is some $a>0$ and $0<b<1$ with $\|\psi^j\|_{\mathcal{L}} \leq ab^j$ for every $j\geq 0$ by \cite{Bosq2000}, Lemma 3.1, and since  
$\nu_p(A(Y)) \leq \|A\|_{\mathcal{L}}\nu_p(Y)$ for any operator $A\in\mathcal{L_H}$ and random variable $Y\in L^p_{\mathcal{H}},$ for any $k\in\mathbb{Z}$ holds
\begin{align*}
\nu_p(X_k)&\leq  a\sum^\infty_{j=0}b^j\nu_p(\varepsilon_{k-j}) \leq  \nu_p(\varepsilon_0)\,\frac{a}{1-b} < \infty.
\end{align*}
In other words, $(X_k)\subset L^p_{\mathcal{H}}.$ Moreover, due to $X^{(m)}_m = \sum^{m-1}_{j=0}\psi^j(\varepsilon_{m-j}) + \sum^{\infty}_{j=m}\psi^j(\varepsilon^{(m)}_{m-j})$ for each $m\in\mathbb{N},$ with $\varepsilon^{(m)}_{m-j}$ being i.i.d. \frenchspacing copies of $\varepsilon_{m-j}$ for each $j,m,$ and because $0< b <1,$ it holds    
\begin{align*}
\sum^\infty_{m=1}\nu_p\big(X_m - X^{(m)}_m\big) 
&\leq  a\sum^\infty_{m=1}\sum^\infty_{j=m}b^j\nu_p\big(\varepsilon_{m-j} - \varepsilon^{(m)}_{m-j}\big)\allowdisplaybreaks\\
&\leq  \nu_p(\varepsilon_0)\,\frac{2a}{1-b}\sum^\infty_{m=1}b^m = \nu_p(\varepsilon_0)\,\frac{2ab}{(1-b)^2}< \infty.
\end{align*}
This shows $L^p$-$m$-approximability of $(X_k).$
\end{proof}

\begin{proof}[\textbf{Proof of Proposition \ref{Prop: Joint stat cart}}] According to Assumption \ref{As: X L^2, centered},  $\boldsymbol{X}=(X_k)$ and $\boldsymbol{Y}=(Y_k)$ are causal processes with finite fourth moments, and as such weakly stationary. Thus, as causality transfers to  $\boldsymbol{X}^{[m]}$ and  $\boldsymbol{Y}^{[n]},$ and   $\|\boldsymbol{z}\|^2=\sum^r_{i=1}\|z_i\|^2$ for $\boldsymbol{z}\coloneqq(z_1, \dots, z_r)^\top\!,$ $r\in\mathbb{N},$ $\boldsymbol{X}^{[m]}$ and  $\boldsymbol{Y}^{[n]}$ are also weakly stationary. Further, from Assumption \ref{As: jointly stationary} follows for each  $h,k,\ell$ and $\boldsymbol{x}=(x_1, \dots, x_m)^\top\!\in\mathcal{H}^m\colon$ 
  \begin{align*}
    \mathscr{C}_{\!X^{[m]}_k, Y^{[n]}_\ell}(x) &= \Erw\!\big\langle X^{[m]}_k, x\big\rangle  Y^{[n]}_\ell\\
    &= \sum^m_{j=1}\Big(\underbrace{\;\Erw\langle X_{k-j+1}, x_j\rangle Y_\ell\;}_{=\,\mathscr{C}_{X_{k+h-j+1}, Y_{\ell+h}}\!(x_j)},\, \underbrace{\;\Erw\langle X_{k-j+1}, x_j\rangle Y_{\ell-1}\;}_{=\,\mathscr{C}_{X_{k+h-j+1}, Y_{\ell+h-1}}\!(x_j)},\, \dots, \underbrace{\;\Erw\langle X_{k-j+1}, x_j\rangle Y_{\ell-n+1}\;}_{=\,\mathscr{C}_{X_{k+h-j+1}, Y_{\ell+h-n+1}}\!(x_j)}\,\Big)^\top\\
    &=\mathscr{C}_{\!X^{[m]}_{k+h}, Y^{[n]}_{\ell+h}}\!(x).
  \end{align*}
Joint weak stationarity of $\boldsymbol{X}^{[m]}$ and  $\boldsymbol{Y}^{[n]}$ is thus verified. Moreover, since $\langle a\otimes b, c\otimes d\rangle_{\mathcal{S}} = \langle a, c\rangle\langle b, d\rangle_\star$ for any $a,c\in\mathcal{H},b,d\in\mathcal{H}_\star,$ due to the definition of $X^{[m]}_s\!, Y^{[n]}_t, \langle\boldsymbol{v},\boldsymbol{w}\rangle$ and $\langle\boldsymbol{x},\boldsymbol{y}\rangle_\star$ for any  $\boldsymbol{v},\boldsymbol{w}\in\mathcal{H}^m$ and $\boldsymbol{x},\boldsymbol{y}\in\mathcal{H}^n_\star,$ and due to Assumption \ref{As: jointly stationary}, for each  $h,i,j,k,\ell, m, n$ holds 
\begin{align*}
  \Erw\!\big\langle X^{[m]}_i\otimes Y^{[n]}_j, X^{[m]}_k\otimes Y^{[n]}_\ell\big\rangle_{\!\mathcal{S}} 
  &= \sum^m_{p=1}\sum^n_{q=1}\,\Erw\langle X_{i-p+1}, X_{k-p+1}\rangle\langle  Y_{j-q+1},  Y_{\ell-q+1}\rangle_\star\allowdisplaybreaks\\
  &= \sum^m_{p=1}\sum^n_{q=1}\;\,\underbrace{\,\Erw\langle X_{i-p+1}\otimes Y_{j-q+1}, X_{k-p+1}\otimes Y_{\ell-q+1}\rangle_{\mathcal{S}}\;}_{=\,\Erw\langle X_{i-p+1+h}\otimes Y_{j-q+1+h},\,X_{k-p+1+h}\otimes Y_{\ell-q+1+h}\rangle_{\mathcal{S}}}\allowdisplaybreaks\\
  &= \Erw\!\big\langle X^{[m]}_{i+h}\otimes Y^{[n]}_{j+h}, X^{[m]}_{k+h}\otimes Y^{[n]}_{\ell+h}\big\rangle_{\!\mathcal{S}}\,.
\end{align*}
Consequently, joint weak stationarity of the fourth moments of $\boldsymbol{X}^{[m]}$ and  $\boldsymbol{Y}^{[n]}$ is shown. 
\end{proof}

\begin{proof}[\textbf{Proof of Theorem \ref{Theo cross-cov op}}] First, let $h\geq 0.$ From the definition of our lag-$h$ cross-covariance operators and their associated empirical versions, Proposition \ref{Prop: Joint stat cart} and elementary conversions  follows with $Z^{[m,n]}_{k,h}\coloneqq X^{[m]}_k\!\otimes\!Y^{[n]}_{k+h} - \mathscr{C}^h_{\!\boldsymbol{X}^{[m]}, \boldsymbol{Y}^{[n]}}$ and $N'=\min\{N,N-h\} - \max\{m,n-h\}+1\colon$
\begin{align}
  N'\Erw\!\big\|\hat{\mathscr{C}}^h_{\!\boldsymbol{X}^{[m]}, \boldsymbol{Y}^{[n]}}\! - \mathscr{C}^h_{\!\boldsymbol{X}^{[m]}, \boldsymbol{Y}^{[n]}}\big\|^2_{\mathcal{S}}\,   &= \frac{1}{N'}\sum^{N'}_{i=1}\sum^{N'}_{j=1}\, \Erw\!\big\langle Z^{[m,n]}_{1,h}, Z^{[m,n]}_{1+k,h} \big\rangle_{\!\mathcal{S}}\notag\allowdisplaybreaks\\
  &= \sum_{|k|<N'}\!\Big(1-\frac{|k|}{N'}\Big)\Erw\!\big\langle Z^{[m,n]}_{1,h}, Z^{[m,n]}_{1+k,h} \big\rangle_{\!\mathcal{S}}\notag\allowdisplaybreaks\\
  &\leq \Erw\!\big\|Z^{[m,n]}_{1,h}\big\|^2_{\mathcal{S}} + 2\sum_{k\geq 1}\,\big|\Erw\!\big\langle Z^{[m,n]}_{1,h}, Z^{[m,n]}_{1+k,h} \big\rangle_{\!\mathcal{S}}\big|.\label{main ineq proof cross-cov}\allowdisplaybreaks
\end{align}
Due to strict stationarity of $(X^{[m]}_k)_k$ and $(Y^{[n]}_k)_k,$ $\|x\otimes y\|_{\mathcal{S}} = \|x\| \|y\|_\star$ for $x\in\mathcal{H}, y\in\mathcal{H}_\star,$ and CSI, 
\begin{align}
  \nu^2_{2,\mathcal{S}}\big(Z^{[m,n]}_{1,h}\big) &= \Erw\!\big\|Z^{[m,n]}_{1,h}\big\|^2_{\mathcal{S}}\, \leq  \big(\Erw\!\big\|X^{[m]}_m\big\|^4\Erw\!\big\|Y^{[n]}_n\big\|^4_\star\big)^{1/2} = \nu^2_4\big(X^{[m]}_m\big)\nu^2_4\big(Y^{[n]}_n\big).\label{4th moment of cart. process, cross-cov}
\end{align}
Further, due to strict stationarity of $(X_k)$ and $(Y_k),$ and because $\|X^{[m]}_k\|^2 = \sum^m_{i=1}\|X_{k-i+1}\|^2$ and $\|Y^{[n]}_k\|^2 = \sum^n_{j=1}\|Y_{k-j+1}\|^2,$ for any $m,n\in\mathbb{N}$ holds 
\begin{align}
    \nu^2_4\big(X^{[m]}_m\big) &\leq m\nu^2_4(X_1)\quad\text{and}\quad\nu^2_4\big(Y^{[n]}_n\big) \leq n\nu^2_4(Y_1).\label{norm of cartesian X,Y, cross-cov}
\end{align}
\noindent To reiterate, after Assumption \ref{As: X L^2, centered} holds $X_k = f(\varepsilon_k, \varepsilon_{k-1}, \dots)$ and $Y_k = g(\eta_k, \eta_{k-1}, \dots)$ for each $k,$ where $f,g$ are some measurable functions, and $(\varepsilon_k)$ and $(\eta_k)$ are some i.i.d. \frenchspacing processes where $\varepsilon_k$ and $\eta_\ell$ are independent for $k\neq \ell.$ Further,  $X^{(\ell)}_k = f(\varepsilon_k, \dots, \varepsilon_{k-\ell+1}, \varepsilon^{(\ell)}_{k-\ell}, \varepsilon^{(\ell)}_{k-\ell-1},   \dots)$ and $Y^{(\ell)}_k = g(\eta_k, \dots, \eta_{k-\ell+1}, \eta^{(\ell)}_{k-\ell}, \eta^{(\ell)}_{k-\ell-1},   \dots)$ for $k\in\mathbb{Z}, \ell\in\mathbb{N}_0,$ where the processes $(\varepsilon^{(\ell)}_k)_k$ and $(\eta^{(\ell)}_k)_k$ are independent copies of $(\varepsilon_k)$ and $(\eta_k),$  respectively. Hereafter, we make the convention that all  $(\varepsilon^{(\ell)}_k)_k$ are independent of $(\eta_k)$ and all $(\eta^{(\ell)}_k)_k,$ and all $(\eta^{(\ell)}_k)_k$ are independent of $(\varepsilon_k)$ and all $(\varepsilon^{(\ell)}_k)_k.$ Further, in \eqref{process X in Cartesian space}--\eqref{approximate of process Y in Cartesian space}, we defined the Cartesian product space-valued processes by $X^{[m]}_k=(X_k, X_{k-1}, \dots, X_{k-m+1})^{\!\top}$ and $Y^{[n]}_k =  (Y_k, Y_{k-1}, \dots, Y_{k-n+1})^{\!\top}\!,$ and its approximates by $  X^{[m], (\ell)}_k= (X^{(\ell)}_k, X^{(\ell)}_{k-1}, \dots, X^{(\ell)}_{k-m+1})^{\!\top}$ and $Y^{[n], (\ell)}_k =(Y^{(\ell)}_k, Y^{(\ell)}_{k-1}, \dots, Y^{(\ell)}_{k-n+1})^{\!\top}\!,$ respectively. A closer look reveals that $Z^{[m,n]}_{1,h}=X^{[m]}_1\!\otimes\!Y^{[n]}_{1+h} - \mathscr{C}^h_{\!\boldsymbol{X}^{[m]},\boldsymbol{Y}^{[n]}}$ and $Z^{[m,n], (\ell)}_{1+k,h}\coloneqq X^{[m],(\ell)}_{1+k}\otimes Y^{[n],(\ell)}_{1+k+h} - \mathscr{C}^h_{\!\boldsymbol{X}^{[m]},\boldsymbol{Y}^{[n]}}$ are  independent for $h\geq 0$ for any $0 \leq \ell\leq 1 + k-h-\max\{m,n-h\}.$ Thus,  \eqref{4th moment of cart. process, cross-cov}--\eqref{norm of cartesian X,Y, cross-cov} and stationarity yield for $\ell'\coloneqq1 + k-h-\max\{m,n-h\}$ and $k'\coloneqq \max\{m,n-h\} + h - 1\colon$ 
\begin{align}
  \big|\Erw\!\big\langle Z^{[m,n]}_{1,h}, Z^{[m,n]}_{1+k,h} \big\rangle_{\!\mathcal{S}}\big| &= \big|\Erw\!\big\langle Z^{[m,n]}_{1,h}, Z^{[m,n]}_{1+k,h} - Z^{[m,n], (\ell')}_{1+k,h}\big\rangle_{\!\mathcal{S}}\big|\notag\allowdisplaybreaks\\
  &\leq \sqrt{mn}\,\nu_4(X_1)\,\nu_4(Y_1)\,\nu_{2,\mathcal{S}}\big(Z^{[m,n]}_{k,h} -  Z^{[m,n],(k)}_{k,h}\big), \quad k\geq \max\{1,k'\}.\label{ineq summands, k>h}
\end{align}
Furthermore, due to $\|a\otimes b - c\otimes d\|^2_{\mathcal{S}} \leq 2 \|a\|^2\|b-d\|^2_\star + 2\|a-c\|^2\|d\|^2_\star$ for $a,c \in\mathcal{H},$ $b, d \in \mathcal{H}_\star,$ strict stationarity,  the CSI, and because of \eqref{norm of cartesian X,Y, cross-cov}, we have
\begin{align}
  \nu_{2,\mathcal{S}}\big(Z^{[m,n]}_{k,h} -  Z^{[m,n],(k)}_{k,h}\big) &\leq  \sqrt{2}\, \nu_4\big(X^{[m]}_m\big)\,\nu_4\big(Y^{[n]}_k\! - Y^{[n],(k)}_k\big) + \sqrt{2}\,\nu_4\big(X^{[m]}_k\! - X^{[m],(k)}_k\big)\,\nu_4\big(Y^{[n]}_n\big)\notag\\
  &\leq \sqrt{2mn}\,\Big[\nu_4(X_1)\,\nu_4\big(Y_k - Y^{(k)}_k\big) + \nu_4(Y_1)\,\nu_4\big(X_k - X^{(k)}_k\big)\Big].\label{ineq tensor Lp-m, cross-cov}
\end{align}
Further, due to the CSI and \eqref{4th moment of cart. process, cross-cov}--\eqref{norm of cartesian X,Y, cross-cov}:
\begin{align}
  \big|\Erw\!\big\langle Z^{[m,n]}_{1,h}, Z^{[m,n]}_{1+k,h} \big\rangle_{\!\mathcal{S}}\big| &\leq \nu_{2,\mathcal{S}}\big(Z^{[m,n]}_{1,h}\big)\,\nu_{2,\mathcal{S}}\big(Z^{[m,n]}_{1+k,h}\big) \leq mn\,\nu^2_4(X_1)\,\nu^2_4(Y_1), \quad 1\leq k<k'.\label{ineq summands, k leq h}
\end{align}

For $h<0,$ \eqref{main ineq proof cross-cov}--\eqref{ineq summands, k leq h} also holds, but since $Z^{[m,n]}_{1,h}$ and $Z^{[m,n], (\ell)}_{1+k,h}$ are then independent  for $0 \leq \ell\leq 1 + k-\max\{m,n-h\},$ both $k'$ and $\ell'$ in \eqref{ineq summands, k>h}--\eqref{ineq tensor Lp-m, cross-cov} have to be replaced by $k''\coloneqq \max\{m,n-h\}- 1$ and $\ell''\coloneqq 1 + k-\max\{m,n-h\},$ respectively. 

Overall, since $\kappa= \max\{m,n-h\} + \mathds 1_{\mathbb{N}}(h) \cdot h - 1$ matches with $k'$ and $k''$ for $h\geq 0$ and $h<0,$ respectively, with $\ell_\kappa\coloneqq 1+k-\max\{m,n-h\} - \mathds 1_{\mathbb{N}}(h) \cdot h - 1$ for any $h\in\mathbb{Z}$ instead of $\ell'$ and $\ell''$ for $h\geq 0$ and $h<0,$ respectively, \eqref{claim cross-cov operator} is considered verified for any $h\in\mathbb{Z}.$ 
%At last, as the upper bound of $\|\hat{\mathscr{C}}^{h}_{\!\boldsymbol{Y}^{[n]}, \boldsymbol{X}^{[m]}}\! - \mathscr{C}^{h}_{\!\boldsymbol{Y}^{[n]}, \boldsymbol{X}^{[m]}}\|_{\mathcal{S}}$ would be the same for $h\geq 0,$ and since $\|\hat{\mathscr{C}}^h_{\!\boldsymbol{X}^{[m]}, \boldsymbol{Y}^{[n]}}\! - \mathscr{C}^h_{\!\boldsymbol{X}^{[m]}, \boldsymbol{Y}^{[n]}}\|_{\mathcal{S}} =  \|\hat{\mathscr{C}}^{-h}_{\!\boldsymbol{Y}^{[n]}, \boldsymbol{X}^{[m]}}\! - \mathscr{C}^{-h}_{\!\boldsymbol{Y}^{[n]}, \boldsymbol{X}^{[m]}}\|_{\mathcal{S}}$ for any $h\in\mathbb{Z},$ \eqref{claim cross-cov operator} also holds for $h < 0.$ 
\end{proof}

\begin{proof}[\textbf{Proof of Proposition \ref{Corollary Cross-cov on X's}}] The proof is almost identical to the proof of Theorem \ref{Theo cross-cov op}, where Assumption \ref{As: jointly stationary} follows from Assumption  \ref{As: X L^2, centered}, since  we are given only the   strictly stationary time series $\boldsymbol{X}=(X_k)$ with finite fourth moments. Here, however, instead of \eqref{ineq tensor Lp-m, cross-cov}, we get for $Z^{[m,n]}_{1,h}=X^{[m]}_1\!\otimes\!X^{[n]}_{1+h} - \mathscr{C}^h_{\!\boldsymbol{X}^{[m]},\boldsymbol{X}^{[n]}}$ and $Z^{[m,n], (\ell)}_{1+k,h}= X^{[m],(\ell)}_{1+k}\otimes X^{[n],(\ell)}_{1+k+h} - \mathscr{C}^h_{\!\boldsymbol{X}^{[m]},\boldsymbol{X}^{[n]}}$  for any $h,k,m,n\colon$
\begin{align*}
  \nu^2_{2,\mathcal{S}}\big(Z^{[m,n]}_{k,h} -  Z^{[m,n],(k)}_{k,h}\big) &\leq  2 \Erw\!\big\|X^{[m]}_{k}\big\|^2\big\|X^{[m]}_{k+h} - X^{[m],(k)}_{k+h}\big\|^2 + 2\Erw\!\big\|X^{[m]}_{k} - X^{[m],(k)}_{k}\big\|^2\big\|X^{[m]}_{k+h}\big\|^2\notag\allowdisplaybreaks\\
  &\leq 4\,\nu^2_4\big(X^{[m]}_m\big)\nu^2_4\big(X^{[m]}_k - X^{[m],(k)}_k\big)\notag\allowdisplaybreaks\\
  &\leq 4\,m^2\nu^2_4(X_1)\nu^2_4\big(X_k - X^{(k)}_k\big).\label{ineq tensor Lp-m}
\end{align*}
By using this inequality and all the other arguments in the proof of Theorem \ref{Theo cross-cov op} with $\boldsymbol{X}=\boldsymbol{Y}$ a.s., the postulated claim is verified.
\end{proof}

\begin{proof}[\textbf{Proof of Corollary \ref{Cor cov op}}] Follows immediately from Proposition \ref{Corollary Cross-cov on X's} with $m=n.$
%We use all inequalities and conversions from the proof of Theorem \ref{Theo cross-cov op} except for  \eqref{ineq tensor Lp-m, cross-cov}, where we put $Z^{[m,n]}_{k,h}=X^{[m]}_k\otimes X^{[n]}_{k+h} - \mathscr{C}^h_{\!\boldsymbol{X}^{[m]}, \boldsymbol{X}^{[n]}},$ and where we replace $\nu_4(Y^{[n]}_n)$ by $\nu_4(X^{[n]}_n),$ and $\nu_4(Y_1)$ by $\nu_4(X_1).$ Further, similar as in the proof of Theorem \ref{Cor cov op}, we obtain due to $\|a\otimes b - c\otimes d\|^2_{\mathcal{S}} \leq 2 \|a\|^2\|b-d\|^2 + 2\|a-c\|^2\|d\|^2$ for $a,b,c,d \in\mathcal{H},$ and by using stationarity, the CSI, and \eqref{4th moment of cart. process, cross-cov}--\eqref{norm of cartesian X,Y, cross-cov}:
%\begin{align*}
%  \nu^2_{2,\mathcal{S}}\big(Z^{[m,n]}_{k,h} -  Z^{[m,n],(k)}_{k,h}\big) &\leq  2\, \nu^2_4\big(X^{[m]}_m\big)\,\nu^2_4\big(X^{[n]}_k\! - X^{[n],(k)}_k\big) + 2\,\nu^2_4\big(X^{[m]}_k\! - X^{[m],(k)}_k\big)\,\nu^2_4\big(X^{[n]}_n\big)\\
%  &\leq 4mn\,\nu_4(X_1)\,\nu_4\big(X_k - X^{(k)}_k\big).
%\end{align*}
%With this inequality, and using all the other inequalities mentioned aboved, the claim is verified. 
\end{proof}

\end{document}